\documentclass[11pt]{article} 
\usepackage{amsmath,amssymb,amsthm}
\usepackage{fancybox}  
\usepackage{graphicx}
\usepackage{color}
\usepackage{colortab}
\usepackage{colortbl}
\usepackage{mathdots}

\allowdisplaybreaks

\begin{document}

\def\fl#1{\left\lfloor#1\right\rfloor}
\def\cl#1{\left\lceil#1\right\rceil}
\def\ang#1{\left\langle#1\right\rangle}
\def\stf#1#2{\left[#1\atop#2\right]} 
\def\sts#1#2{\left\{#1\atop#2\right\}}
\def\eul#1#2{\left\langle#1\atop#2\right\rangle}
\def\N{\mathbb N}
\def\Z{\mathbb Z}
\def\R{\mathbb R}
\def\C{\mathbb C}
\newcommand{\ctext}[1]{\raise0.2ex\hbox{\textcircled{\scriptsize{#1}}}}

\newtheorem{theorem}{Theorem}
\newtheorem{Prop}{Proposition}
\newtheorem{Cor}{Corollary}
\newtheorem{Lem}{Lemma}
\newtheorem{Sublem}{Sublemma}
\newtheorem{Def}{Definition}  
\newtheorem{Conj}{Conjecture}

\newenvironment{Rem}{\begin{trivlist} \item[\hskip \labelsep{\it
Remark.}]\setlength{\parindent}{0pt}}{\end{trivlist}}

\title{$p$-numerical semigroups of Pell triples}

\author{
Takao Komatsu 
\\
\small Department of Mathematical Sciences, School of Science\\[-0.8ex]
\small Zhejiang Sci-Tech University\\[-0.8ex]
\small Hangzhou 310018 China\\[-0.8ex]
\small \texttt{komatsu@zstu.edu.cn}\\\\
Jiaxin Mu\\
\small Department of Mathematical Sciences, School of Science\\[-0.8ex]
\small Zhejiang Sci-Tech University\\[-0.8ex]
\small Hangzhou 310018 China\\[-0.8ex] 
\small \texttt{jiaxinm68@163.com}  
}

\date{
%\small Submitted: February 15, 2023;  Accepted: March 25, 2023.\\
\small MR Subject Classifications: Primary 11D07; Secondary 20M14, 05A17, 05A19, 11D04, 11B68, 11P81 
}

\maketitle
 
\begin{abstract} 
For a nonnegative integer $p$, the $p$-numerical semigroup $S_p$ is defined as the set of integers whose nonnegative integral linear combinations of given positive integers $a_1,a_2,\dots,a_\kappa$ with $\gcd(a_1,a_2,\dots,a_\kappa)=1$ are expressed in more than $p$ ways. When $p=0$, $S=S_0$ is the original numerical semigroup. The largest element and the cardinality of $\mathbb N_0\backslash S_p$ are called the $p$-Frobenius number and the $p$-genus, respectively. Their explicit formulas are known for $\kappa=2$, but those for $\kappa\ge 3$ have been found only in some special cases. For some known cases, such as the Fibonacci and the Jacobsthal triplets, similar techniques could be applied and explicit formulas such as the $p$-Frobenius number could be found.  

In this paper, we give explicit formulas for the $p$-Frobenius number and the $p$-genus of Pell numerical semigroups $\bigl(P_i(u),P_{i+2}(u),P_{i+k}(u)\bigr)$. Here, for a given positive integer $u$, Pell-type numbers $P_n(u)$ satisfy the recurrence relation $P_n(u)=u P_{n-1}(u)+P_{n-2}(u)$ ($n\ge 2$) with $P_0(u)=0$ and $P_1(u)=1$. The $p$-Ap\'ery set is used to find the formulas, but it shows a different pattern from those in the known results, and some case by case discussions are necessary.   
\end{abstract}

\section{Introduction}  

Given the set of positive integers $A:=\{a_1,a_2,\dots,a_\kappa\}$ ($\kappa\ge 2$), for a nonnegative integer $p$, let $S_p$ be the set of integers whose nonnegative integral linear combinations of given positive integers $a_1,a_2,\dots,a_\kappa$ are expressed in more than $p$ ways.  
For some backgrounds of the number of representations, see, e.g., \cite{Binner20,Cayley,Ko03,sy1857,tr00}.  
For a set of nonnegative integers $\mathbb N_0$, the set $\mathbb N_0\backslash S_p$ is finite if and only if $\gcd(a_1,a_2,\dots,a_\kappa)=1$.  Then there exists the largest integer $g_p(A):=g(S_p)$ in $\mathbb N_0\backslash S_p$, which is called the {\it $p$-Frobenius number}. The cardinality of $\mathbb N_0\backslash S_p$ is called the {\it $p$-genus} and is denoted by $n_p(A):=n(S_p)$.   
The sum of the elements in $\mathbb N_0\backslash S_p$ is called the {\it $p$-Sylvester sum} and is denoted by $s_p(A):=s(S_p)$.  
 This kind of concept is a generalization of the famous Diophantine problem of Frobenius (\cite{RA05,sy1882,sy1884}) since $p=0$ is the case when the original Frobenius number $g(A)=g_0(A)$, the genus $n(A)=n_0(A)$   
and the Sylvester sum $s(A)=s_0(A)$   
are recovered. We can call $S_p$ the {\it $p$-numerical semigroup}. Strictly speaking, when $p\ge 1$, $S_p$ does not include $0$ since the integer $0$ has only one representation, so it satisfies simply additivity and the set $S_p\cup\{0\}$ becomes a numerical semigroup. For numerical semigroups, we refer to \cite{ADG,RG1,RG2}. For the $p$-numerical semigroup, we refer to \cite{KY24}.   

There exist different extensions of the Frobenius number and genus, even in terms of the number of representations called denumerant. For example, some considered $S_p^\ast$ as the set of integers whose nonnegative integral linear combinations of given positive integers $a_1,a_2,\dots,a_\kappa$ are expressed in {\it exactly} $p$ ways (see, e.g., \cite{BDFHKMRSS,FS11}). Then, 
the corresponding $p$-Frobenius number $g_p^\ast(A)$ is the largest integer that has {\it\underline{exactly $p$ distinct}} representations. However, in this case, $g_p^\ast(A)$ is not necessarily increasing as $p$ increases. For example, when $A:=\{2,5,7\}$, $g_{17}^\ast(2,5,7)=43>g_{18}^\ast(2,5,7)=42$. In addition, for some $p$, $g_p^\ast$ may not exist. For example, $g_{22}^\ast(2,5,7)$ does not exist because there is no positive integer whose number of representations is exactly $22$.  
The $p$-genus may be also defined in different ways. For example, $n_p^\ast(A)$ can be defined the cardinality of $[\ell_p(A),g_p(A)+1]\backslash S_p(A)$, where $\ell_p(A)$ is the least element of $S_p(A)$. However, in our definition of $n_p(A)$ as the cardinality of   $[0,g_p(A)+1]\backslash S_p(A)$, we can use the convenient formula arising from the $p$-Ap\'ery set in order to obtain $n_p(A)$. See the next section.  
      
%Consider the three-term recurrence relation $\mathcal W_n(u,v)=u\mathcal W_{n-1}(u,v)+v\mathcal W_{n-2}(u,v)$ ($n\ge 2$) with $\mathcal W_0(u,v)=0$ and $\mathcal W_1(u,v)=1$. When $u=v=1$, $F_n=\mathcal W_n(1,1)$ are Fibonacci numbers. When $u=1$ and $v=2$, $J_n=\mathcal W_n(1,2)$ are Jacobsthal numbers. When $u=2$ and $v=1$, $P_n=\mathcal W_n(1,1)$ are Pell numbers. 
In \cite{Matt09}, numerical semigroups generated by $\{a,a+b,a F_{k-1}+b F_k\}$ are considered. Using a technique of Johnson \cite{Johnson60}, the Frobenius numbers of such semigroups are found as a generalization of the result on the set $\{F_i,F_{i+2},F_{i+k}\}$ by Marin et. al. \cite{MAR}. Some more generalizations using Fibonacci numbers have been done in \cite{SKT,Fel09}. The results of replacing the Fibonacci number with the Luca number (with the same recurrence and only different initial values) have been found fairly similarly (\cite{BYM20}), but replacing the Fibonacci number with a generalized number have been very difficult. 
However, it has only recently succeeded in obtaining the result of replacing the Fibonacci number with the generalized Jacobsthal number (Jacobsthal polynomial). This includes the result for the Fibonacci number. In addition, the results in the Jacobsthal-Lucas polynomial including Lucas numbers have been similarly determined.    
However, contrary to easy predictions and wishful observations by many people, it has been more difficult to find the results for the Pell number. One of the reasons for this is that the Pell number grows more faster than the Jacobsthal number, and it is more difficult to grasp the actual situation. In addition, case by case discussions are more delicate, and it is difficult to obtain general properties, which further makes handling difficult.   

In this paper, for a fixed positive integer $u$, we treat with Pell-type numbers (Pell polynomials) $P_n(u)$, defined by the recurrence relation $P_n(u)=u P_{n-1}(u)+P_{n-2}(u)$ ($n\ge 2$) with $P_0(u)=0$ and $P_1(u)=1$ (see, e.g., \cite[Chapter 31]{Koshy2}). When $u=1$, $F_n=P_n(1)$ are Fibonacci numbers. When $u=2$, $P_n=P_n(2)$ are the classical Pell numbers. Then, 
we give explicit formulas of $p$-Frobenius numbers for $A:=\{P_i(u),P_{i+2}(u),P_{i+k}(u)\}$ ($i,k\ge 3$).

For $\kappa=2$, that is, the case of two variable, closed formulas are explicitly given for $g_0(A)$ (\cite{sy1884}), $n_0(A)$ (\cite{sy1882}) and $s_0(A)$ (\cite{bs93}; its extension \cite{ro94}). 
However, for $\kappa\ge 3$, the Frobenius number cannot be given by any set of closed formulas which can be reduced to a finite set of certain polynomials (\cite{cu90}). For $k=3$, various algorithms have been devised for finding the Frobenius number (\cite{Fel06,Johnson60,RGS04}). Some inexplicit formulas for the Frobenius number in three variables can be seen in \cite{tr17}.   
The need to care enough is that even in the original case of $p=0$, it is very difficult to give a closed explicit formula of any general sequence for three or more variables (see, e.g., \cite{RR18,RGS04,RBT2015,RBT2017}). Indeed, it is even more difficult when $p>0$. However, finally, we have succeeded in obtaining $p$-Frobenius number in triangular numbers \cite{Ko22a} and repunits \cite{Ko22b} as well as Fibonacci and Lucas triplets \cite{KY}, Jacobsthal triples \cite{KP,KLP} and arithmetic triples \cite{KY23c} quite recently. 

In the case of the Pell number, it was thought that the formula would be easily found in a similar manner to these in \cite{Ko22b,KY,KP,KLP}, but things were not simple. The reason for this is that the Pell number increased faster than the Fibonacci and Jacobsthal numbers, making it difficult to predict the results experimentally. 
In addition, it is also hinted that the combination of the triplets was likely to show completely different behavior.   
However, in this paper it is finally successful to find the $p$-Frobenius number and the $p$-genus for the more general Pell polynomial triplets. The situation is different in each case where the triplets are (even,even,odd), (odd,odd,even), and (odd,odd,odd). In particular, it is noted that in the case of (even,even,odd) and (odd,odd,even), we find new processes that have not been seen in previously published papers.

\section{Preliminaries}  

We introduce the Ap\'ery set (see \cite{Apery}) below in order to obtain the formulas for $g_p(A)$, $n_p(A)$ and $s_p(A)$ technically. Without loss of generality, we assume that $a_1=\min(A)$. 

\begin{Def}  
Let $p$ be a nonnegative integer. For a set of positive integers $A=\{a_1,a_2,\dots,a_\kappa\}$ with $\gcd(A)=1$ and $a_1=\min(A)$ we denote by 
$$
{\rm Ap}_p(A)={\rm Ap}_p(a_1,a_2,\dots,a_\kappa)=\{m_0^{(p)},m_1^{(p)},\dots,m_{a_1-1}^{(p)}\}\,, 
$$  
the $p$-Ap\'ery set of $A$, where each positive integer $m_i^{(p)}$ $(0\le i\le a_1-1)$ satisfies the conditions:  
$$
{\rm (i)}\, m_i^{(p)}\equiv i\pmod{a_1},\quad{\rm (ii)}\, m_i^{(p)}\in S_p(A),\quad{\rm (iii)}\, m_i^{(p)}-a_1\not\in S_p(A)\,. 
$$ 
Note that $m_0^{(0)}$ is defined to be $0$.  
\label{apery} 
\end{Def}  

\noindent 
It follows that for each $p$, 
$$
{\rm Ap}_p(A)\equiv\{0,1,\dots,a_1-1\}\pmod{a_1}\,. 
$$  

Even though it is hard to find any explicit form of $g_p(A)$ as well as $n_p(A)$ and $s_p(A)$ $k\ge 3$, by using convenient formulas established in \cite{Ko22,Ko-p}, we can obtain such values for some special sequences $(a_1,a_2,\dots,a_\kappa)$ after finding any regular structure of $m_j^{(p)}$. One convenient formula is on the power sum 
$$
s_p^{(\mu)}(A):=\sum_{n\in\mathbb N_0\backslash S_p(A)}n^\mu
$$  
by using Bernoulli numbers $B_n$ defined by the generating function 
$$
\frac{x}{e^x-1}=\sum_{n=0}^\infty B_n\frac{x^n}{n!}\,, 
$$
and another convenient formula is on the weighted power sum (\cite{KZ0,KZ}) 
$$
s_{\lambda,p}^{(\mu)}(A):=\sum_{n\in \mathbb N_0\backslash S_p(A)}\lambda^n n^\mu 
$$  
by using Eulerian numbers $\eul{n}{m}$ appearing in the generating function 
$$ 
\sum_{k=0}^\infty k^n x^k=\frac{1}{(1-x)^{n+1}}\sum_{m=0}^{n-1}\eul{n}{m}x^{m+1}\quad(n\ge 1)
$$ 
with $0^0=1$ and $\eul{0}{0}=1$. Here, $\mu$ is a nonnegative integer and $\lambda\ne 1$. From these convenient formulas, many useful expressions are yielded as special cases. 
Some useful ones are given as follows.  The formulas (\ref{mp-n}) and (\ref{mp-s}) are entailed from $s_{\lambda,p}^{(0)}(A)$ and $s_{\lambda,p}^{(1)}(A)$, respectively.  

\begin{Lem}  
Let $\kappa$, $p$ and $\mu$ be integers with $\kappa\ge 2$ and $p\ge 0$.  
Assume that $\gcd(a_1,a_2,\dots,a_\kappa)=1$.  We have 
\begin{align}  
g_p(a_1,a_2,\dots,a_\kappa)&=\left(\max_{0\le j\le a_1-1}m_j^{(p)}\right)-a_1\,, 
\label{mp-g}
\\  
n_p(a_1,a_2,\dots,a_\kappa)&=\frac{1}{a_1}\sum_{j=0}^{a_1-1}m_j^{(p)}-\frac{a_1-1}{2}\,,  
\label{mp-n}
\\
s_p(a_1,a_2,\dots,a_\kappa)&=\frac{1}{2 a_1}\sum_{j=0}^{a_1-1}\bigl(m_j^{(p)}\bigr)^2-\frac{1}{2}\sum_{j=0}^{a_1-1}m_j^{(p)}+\frac{a_1^2-1}{12}\,.
\label{mp-s}
\end{align} 
\label{lem-mp}
\end{Lem} 

\noindent 
{\it Remark.}  
When $p=0$, the formulas (\ref{mp-g}), (\ref{mp-n}) and (\ref{mp-s}) reduce to the formulas by  
Brauer and Shockley \cite[Lemma 3]{bs62}, Selmer \cite[Theorem]{se77}, and Tripathi \cite[Lemma 1]{tr08}\footnote{There was a typo, but it was corrected in \cite{PT}.}, respectively: 
\begin{align*}   
g(a_1,a_2,\dots,a_\kappa)&=\left(\max_{0\le j\le a_1-1}m_j\right)-a_1\,,\\  
n(a_1,a_2,\dots,a_\kappa)&=\frac{1}{a_1}\sum_{j=0}^{a_1-1}m_j-\frac{a_1-1}{2}\,,\\  
s(a_1,a_2,\dots,a_\kappa)&=\frac{1}{2 a_1}\sum_{j=0}^{a_1-1}(m_j)^2-\frac{1}{2}\sum_{j=0}^{a_1-1}m_j+\frac{a_1^2-1}{12}\,,   
\end{align*} 
where $m_j=m_j^{(0)}$ ($1\le j\le a_1-1$) with $m_0=m_0^{(0)}=0$.

\section{Main results}  

We use the trivial formula repeatedly: 
\begin{equation}
P_{i+k}=P_{i+1}P_{k}+P_{i}P_{k-1}\,.  
\label{pell-equ}
\end{equation}

\subsection{Initial results for even $i$}  

Note that (even,even,odd) case does not appear when $u=1$. Let $u\ge 2$.  

First, from the recurrence relation $P_{n}(u)=u P_{n-1}(u)+P_{n-2}(u)$ ($n\ge 2$) with $P_0(u)=0$ and $P_1(u)=1$, when $u\ge 2$, 
\begin{equation}
P_{2 n}(u)\equiv 0\pmod{u}\quad\hbox{and}\quad P_{2 n+1}(u)\equiv 1\pmod{u}\,. 
\label{pell-con-res}
\end{equation} 
Hence, concerning the set $A:=\{P_i(u),P_{i+2}(u),P_{i+k}(u)\}$, if $i$ is even, then $k$ must be odd. Otherwise, $\gcd(A)\ne 1$.  
Since $P_2(u)|P_4(u)$, we consider the triples $A:=\{P_i(u),P_{i+2}(u),P_{i+k}(u)\}$ for $i\ge 3$.  

For integers $i$ and $k$ with $2\le i\le k$, we have\footnote{When $i>k$, the case by case discussion is necessary as $k=i-1$, $k=i-2$ and so on. As it is complicated, we do not treat with the case when $i>k$.}   
\begin{multline*}
g\bigl(P_{2 i}(u),P_{2 i+2}(u),P_{2i+2 k+1}(u)\bigr)\\
=\left(\frac{P_{2 i}(u)}{u}-1\right)P_{2 i+2}(u)+(u-1)P_{2 i+2 k+1}(u)-P_{2 i}(u)\,. 
\end{multline*}  
Note that $P_{2 i}(u)/u$ is a positive integer. In order to prove this result, we must show that all the elements of ${\rm Ap}(A)$ are distributed exactly inside of the rectangle in Table \ref{tb:ap0}.

\begin{table}[htbp]
  \centering
%\scalebox{0.7}{
\begin{tabular}{ccccc}
\cline{1-2}\cline{3-4}\cline{5-5}
\multicolumn{1}{|c}{$t_{0,0}$}&$t_{1,0}$&$\cdots$&$\cdots$&\multicolumn{1}{c|}{$t_{P_{2 i}(u)/u-1,0}$}\\
\multicolumn{1}{|c}{$t_{0,1}$}&$t_{1,1}$&$\cdots$&$\cdots$&\multicolumn{1}{c|}{$t_{P_{2 i}(u)/u-1,1}$}\\
\multicolumn{1}{|c}{$\vdots$}&$\vdots$&$\cdots$&$\cdots$&\multicolumn{1}{c|}{$\vdots$}\\
\multicolumn{1}{|c}{$t_{0,u-1}$}&$t_{1,u-1}$&$\cdots$&$\cdots$&\multicolumn{1}{c|}{$t_{P_{2 i}(u)/u-1,u-1}$}\\
\cline{1-2}\cline{3-4}\cline{5-5}
\end{tabular}
%} 
  \caption{${\rm Ap}_0\bigl(P_{2 i}(u),P_{2 i+2}(u),P_{2i+2 k+1}(u)\bigr)$}
  \label{tb:ap0}
\end{table}

Here, we consider the expression 
$$
t_{y,z}:=y P_{2 i+2}(u)+z P_{2i+2 k+1}(u)\,. 
$$
First, the number of the elements in the rectangle is 
$$
\frac{P_{2 i}(u)}{u}\times u=P_{2 i}(u)\,,
$$
which is equal to the least element of ${\rm Ap}(A)$. 
Next, by $P_{2 i+2}(u)\equiv u P_{2 i-1}(u)\pmod{P_{2 i}(u)}$ %and $P_{2 i+2 k+1}\equiv P_{2 i-1}P_{2 k+1}$
$$
t_{\ell,z}\equiv u\ell P_{2 i-1}(u)+z P_{2i-1}P_{2k+1}\pmod{P_{2 i}(u)}\,. 
$$
If for integers $\ell$ and $\ell'$ with $0\le \ell,\ell'<P_{2 i}(u)/u$  
$$ 
t_{\ell,z}\equiv t_{\ell',z}\pmod{P_{2 i}(u)}\,, 
$$ 
by 
\begin{equation} 
\gcd\bigl(P_{2 i-1}(u),P_{2 i}(u)/u\bigr)=1
\label{eq:p123} 
\end{equation}  
we have 
$$
(\ell-\ell')P_{2 i-1}(u)\equiv 0\left({\rm mod\,}\frac{P_{2 i}(u)}{u}\right)\,. 
$$ 
It entails that  
$\ell\equiv\ell'\left({\rm mod\,}\frac{P_{2 i}(u)}{u}\right)$, that is $\ell=\ell'$.  
Hence, for a fixed $z$ with $0\le z\le u-1$, all the elements of the form $t_{\ell,z}$ are distinct among the complete residue system $\{0,1,\dots,P_{2 i}(u)-1\}$. 
Concerning the fact in (\ref{eq:p123}), if $\gcd\bigl(P_{2 i-1}(u),P_{2 i}(u)/u\bigr)=d$, then 
$d|P_{2 i}(u)/u$ and $d|P_{2 i-1}(u)$, so by the recurrence relation we get 
$d|P_{2 i-2}(u)/u$ and $d|P_{2 i-3}(u)$, finally $d|P_{2}(u)/u=1$ and $d|P_1(u)=1$. Therefore, $d=1$. 

By the condition (\ref{pell-con-res}),  
we have for fixed $z$ ($0\le z\le u-1$) 
$$  
t_{y,z}=y P_{2 i+2}(u)+z P_{2 i+2 k+1}(u)\equiv y\cdot 0+z\cdot 1=z\pmod u\,. 
$$ 
Since $u|P_{2 i}(u)$, any of two elements belonging to the different rows are distinct modulo $P_{2 i}(u)$.  
Hence, we can conclude that all the elements in the rectangle in Table \ref{tb:ap0} constitutes the complete residue system $\{0,1,\dots,P_{2 i}(u)-1\}$, thus the Ap\'ery set ${\rm Ap}(A)$. 

Since the largest element of ${\rm Ap}(A)$ is $t_{\frac{P_{2 i}(u)}{u}-1,u-1}$, by the formula in Lemma \ref{lem-mp} (\ref{mp-g}), we have 
\begin{multline*}
g\bigl(P_{2 i}(u),P_{2 i+2}(u),P_{2i+2 k+1}(u)\bigr)\\
=\left(\frac{P_{2 i}(u)}{u}-1\right)P_{2 i+2}(u)+(u-1)P_{2 i+2 k+1}(u)-P_{2 i}(u)\,. 
\end{multline*}

\subsection{$p=1$}  

When $p=1$, all elements of $1$-Ap\'ery set appear in the form of displacing the elements of $0$-Ap\'ery set  by shift as they are, shown in Table \ref{tb:ap0-1}.  

\begin{table}[htbp]
  \centering
\scalebox{0.7}{                  
  \begin{tabular}{cccccccccccc}
\cline{1-2}\cline{3-4}\cline{5-6}\cline{7-8}\cline{9-10}\cline{11-12}
\multicolumn{1}{|c}{$t_{0,0}$}&$t_{1,0}$&$\cdots$&&$\cdots$&\multicolumn{1}{c|}{$t_{P_{2 i}(u)/u-1,0}$}&$t_{P_{2 i}(u)/u,0}$&$t_{P_{2 i}(u)/u+1,0}$&$\cdots$&&$\cdots$&\multicolumn{1}{c|}{$t_{2 P_{2 i}(u)/u-1,0}$}\\
\multicolumn{1}{|c}{$t_{0,1}$}&$t_{1,1}$&$\cdots$&&$\cdots$&\multicolumn{1}{c|}{$t_{P_{2 i}(u)/u-1,1}$}&$t_{P_{2 i}(u)/u,1}$&$t_{P_{2 i}(u)/u+1,1}$&$\cdots$&&$\cdots$&\multicolumn{1}{c|}{$t_{2 P_{2 i}(u)/u-1,1}$}\\
\multicolumn{1}{|c}{$\vdots$}&$\vdots$&&&&\multicolumn{1}{c|}{$\vdots$}&$\vdots$&$\vdots$&&&&\multicolumn{1}{c|}{$\vdots$}\\
\multicolumn{1}{|c}{$t_{0,u-1}$}&$t_{1,u-1}$&$\cdots$&&$\cdots$&\multicolumn{1}{c|}{$t_{P_{2 i}(u)/u-1,u-1}$}&$t_{2 P_{2 i}(u)/u,u-1}$&$t_{P_{2 i}(u)/u+1,u-1}$&$\cdots$&&$\cdots$&\multicolumn{1}{c|}{$t_{P_{2 i}(u)/u-1,u-1}$}\\
\cline{1-2}\cline{3-4}\cline{5-6}\cline{7-8}\cline{9-10}\cline{11-12}
  \end{tabular}
} 
  \caption{${\rm Ap}_0(A)$ and ${\rm Ap}_1(A)$}
  \label{tb:ap0-1}
\end{table}

Since 
$$
t_{y,z}\equiv t_{y+P_{2 i}(u)/u,z}\pmod{P_{2 i}(u)}\quad (0\le y\le P_{2 i}(u)/u-1,\,0\le z\le u-1)\,,  
$$ 
all these elements constitute a complete residue system modulo by $P_{2 i}(u)$. 
Let the whole expression by 
$$
t_{x,y,z}:=x P_{2 i}(u)+y P_{2 i+2}(u)+z P_{2i+2 k+1}(u)\,. 
$$     
Or, for short, $(x,y,z)$.   
Then, each element is expressed exactly in two ways: 
$$
\left(0,\frac{P_{2 i}(u)}{u}+j,h\right)\quad\hbox{and}\quad \left(\frac{P_{2 i+2}(u)}{u},j,h\right)\,,
$$ 
where $0\le j\le P_{2 i}(u)/u-1$ and $0\le h\le u-1$.  
Therefore, the elements in the area immediately beside the area of ${\rm Ap}_0(A)$ all belong to ${\rm Ap}_1(A)$. It is clear that the largest element of ${\rm Ap}_1(A)$ is $t_{P_{2 i}(u)-1,u-1}$. Hence, by Lemma \ref{lem-mp} (\ref{mp-g}), we have 
\begin{multline*}
g_1\bigl(P_{2 i}(u),P_{2 i+2}(u),P_{2i+2 k+1}(u)\bigr)\\
=\left(2 P_{2 i}(u)/u-1\right)P_{2 i+2}(u)+(u-1)P_{2 i+2 k+1}(u)-P_{2 i}(u)\,. 
\end{multline*}

\subsection{$p\ge 2$}  

When $p=2$, all elements of $2$-Ap\'ery set appear in the form of displacing the elements of $1$-Ap\'ery set by shift. Similarly, all elements of $3$-Ap\'ery set appear in the form of displacing the elements of $2$-Ap\'ery set by shift. Table \ref{tb:g3system400} shows the areas of $p$-Ap\'ery sets, indicated by $\ctext{p}$ for $p=0,1,2,3$.

\begin{table}[htbp]
  \centering
\scalebox{0.7}{
\begin{tabular}{cccccccccccccccc}
\multicolumn{1}{|c}{}&&&\multicolumn{1}{c|}{}&&&&\multicolumn{1}{c|}{}&&&&\multicolumn{1}{c|}{}&&&&\multicolumn{1}{c|}{}\\ 
\cline{1-2}\cline{3-4}\cline{5-6}\cline{7-8}\cline{9-10}\cline{11-12}\cline{13-14}\cline{15-16}
\multicolumn{1}{|c}{}&&&\multicolumn{1}{c|}{}&&&&\multicolumn{1}{c|}{}&&&&\multicolumn{1}{c|}{}&&&&\multicolumn{1}{c|}{}\\ 
\multicolumn{1}{|c}{}&$\ctext{0}$&&\multicolumn{1}{c|}{}&&$\ctext{1}$&&\multicolumn{1}{c|}{}&&$\ctext{2}$&&\multicolumn{1}{c|}{}&&$\ctext{3}$&&\multicolumn{1}{c|}{}\\ 
\cline{1-2}\cline{3-4}\cline{5-6}\cline{7-8}\cline{9-10}\cline{11-12}\cline{13-14}\cline{15-16}
\multicolumn{1}{|c}{}&&&\multicolumn{1}{c|}{}&&&&\multicolumn{1}{c|}{}&&&&\multicolumn{1}{c|}{}&&&&\multicolumn{1}{c|}{}
\end{tabular}
} 
  \caption{${\rm Ap}_p(A)$ ($p=0,1,2,3$)}
  \label{tb:g3system400}
\end{table}

Since for $l=1,2,\dots$ 
$$
t_{y,z}\equiv t_{y+l P_{2 i}(u)/u,z}\pmod{P_{2 i}(u)}\quad (0\le y\le P_{2 i}(u)/u-1,\,0\le z\le u-1)\,,  
$$ 
all these elements constitute a complete residue system modulo by $P_{2 i}(u)$. 

Each element in ${\rm Ap}_p(A)$ has exactly in $p+1$ representations:  
\begin{multline}  
\left(0,\frac{p P_{2 i}(u)}{u}+j,h\right)=\left(\frac{P_{2 i+2}(u)}{u},\frac{(p-1)P_{2 i}(u)}{2}+j,h\right)\\
=\left(P_{2 i+2}(u),\frac{(p-2)P_{2 i}(u)}{u}+j,h\right)= 
\cdots=\left(\frac{p P_{2 i+2}(u)}{u},j,h\right)
\label{eq:ap-rep}
\end{multline} 
where $0\le j\le P_{2 i}(u)/u-1$ and $0\le h\le u-1$.  

Since the largest element of ${\rm Ap}_p(A)$ is $t_{p P_{2 i}(u)/u-1,u-1}$. Hence, by Lemma \ref{lem-mp} (\ref{mp-g}), we have 
\begin{multline*}  
g_p\bigl(P_{2 i}(u),P_{2 i+2}(u),P_{2i+2 k+1}(u)\bigr)\\
=\left(\frac{(p+1)P_{2 i}(u)}{u}-1\right)P_{2 i+2}(u)+(u-1)P_{2 i+2 k+1}(u)-P_{2 i}(u)\,. 
\end{multline*}

However, this situation with regularity is not continued if 
$$
p\ge\frac{u(P_{2 k+1}(u)+1)}{P_{2 i}(u)}-1\,.  
$$  
In this case, in addition to the $p+1$ representations indicated by (\ref{eq:ap-rep}), the value at 
$$
\left(0,\frac{p P_{2 i}(u)}{u}+j,0\right)\quad\hbox{or}\quad \left(0,\frac{p P_{2 i}(u)}{u}+j,1\right)
$$
may have another representation where $j=2$ or $j=3$, respectively. Thus, it can have $p+2$ or more expressions. This implies that a different integer exists, which has the same residue congruent modulo $P_{2 i}(u)$, with $p+1$ representations, and this different integer becomes an element of the $p$-Ap\'ery set. Thus, the Frobenius number cannot keep the regularity after this. 
    
In order to verify this fact, assume that the value at $(x,y,3)$ is equal to that at 
$$
\left(0,\frac{p P_{2 i}(u)}{u}+\frac{P_{2 i}(u)}{u}-1,1\right)\,, 
$$ 
which may become a candidate to take the largest element of ${\rm Ap}_p(A)$. Since by (\ref{pell-equ})
$$
P_{2 i+2 k+1}(u)=\frac{P_{2 i+2}(u)P_{2 k+1}(u)}{u}-\frac{P_{2 i}(u)P_{2 k-1}(u)}{u}\,,
$$ 
such an equality is explained by 
\begin{align*}
&\left(\frac{(p+1)P_{2 i}(u)}{u}-1\right)P_{2 i+2}(u)\\
&=x P_{2 i}(u)+y P_{2 i+2}(u)+P_{2 i+2}(u)P_{2 k+1}(u)-P_{2 i}(u)P_{2 k-1}(u)\\
&=(x-P_{2 k-1}(u))P_{2 i}(u)+(y+P_{2 k+1}(u))P_{2 i+2}(u)\,. 
\end{align*}
Taking $x=P_{2 k-1}(u)$, we have 
$$
y=\frac{(p+1)P_{2 i}(u)}{u}-1-P_{2 k+1}(u)\,. 
$$ 
What there is no additional representation implies that $y<0$, yielding 
$$
p<\frac{u(P_{2 k+1}(u)+1)}{P_{2 i}(u)}-1\,. 
$$

\begin{theorem}
For $2\le i\le k$ and $0\le p<u(P_{2 k+1}(u)+1)/P_{2 i}(u)-1$, we have 
\begin{multline*} 
g_p\bigl(P_{2 i}(u),P_{2 i+2}(u),P_{2i+2 k+1}(u)\bigr)\\
=\left(\frac{(p+1)P_{2 i}(u)}{u}-1\right)P_{2 i+2}(u)+(u-1)P_{2 i+2 k+1}(u)-P_{2 i}(u)\,. 
\end{multline*}   
\label{th1}
\end{theorem}

\subsection{$p$-genus}  

By using Table \ref{tb:ap-p} and Lemma \ref{lem-mp} (\ref{mp-n}) and (\ref{mp-s}), we can obtain the $p$-genus and the $p$-Sylvester sum, respectively, as well.   

\begin{table}[htbp]
  \centering
%\scalebox{0.7}{
\begin{tabular}{ccccc}
\cline{1-2}\cline{3-4}\cline{5-5}
\multicolumn{1}{|c}{$t_{p P_{2 i}(u)/u,0}$}&$t_{p P_{2 i}(u)/u+1,0}$&$\cdots$&$\cdots$&\multicolumn{1}{c|}{$t_{(p+1)P_{2 i}(u)/u-1,0}$}\\
\multicolumn{1}{|c}{$t_{p P_{2 i}(u)/u,1}$}&$t_{p P_{2 i}(u)/u+1,1}$&$\cdots$&$\cdots$&\multicolumn{1}{c|}{$t_{(p+1)P_{2 i}(u)/u-1,1}$}\\
\multicolumn{1}{|c}{$\vdots$}&$\vdots$&&&\multicolumn{1}{c|}{$\vdots$}\\
\multicolumn{1}{|c}{$t_{p P_{2 i}(u)/u,u-1}$}&$t_{p P_{2 i}(u)/u+1,u-1}$&$\cdots$&$\cdots$&\multicolumn{1}{c|}{$t_{(p+1)P_{2 i}(u)/u-1,u-1}$}\\
\cline{1-2}\cline{3-4}\cline{5-5}
\end{tabular}
%} 
  \caption{${\rm Ap}_p\bigl(P_{2 i}(u),P_{2 i+2}(u),P_{2i+2 k+1}(u)\bigr)$}
  \label{tb:ap-p}
\end{table} 

By Table \ref{tb:ap-p}, we have 
\begin{align*}
&\sum_{w\in{\rm Ap}_p(A)}w\\
&=\sum_{y=0}^{P_{2 i}(u)/u-1}\sum_{z=0}^{u-1}\left(\left(\frac{p P_{2 i}(u)}{u}\frac{P_{2 i}(u)}{u}+y\right)u P_{2 i+2}(u)+z\frac{P_{2 i}(u)}{u}P_{2 i+2 k+1}(u)\right)\\
&=\frac{P_{2 i}(u)}{2 u}\left((2 p+1)P_{2 i}(u)P_{2 i+2}(u)+u\bigl(-P_{2 i+2}(u)+(u-1)P_{2 i+2 k+1}(u)\bigr)\right)\,.  
\end{align*} 
Hence, by Lemma \ref{lem-mp} (\ref{mp-n}), we obtain 
\begin{align*}
&n_p\bigl(P_{2 i}(u),P_{2 i+2}(u),P_{2i+2 k+1}(u)\bigr)\\
&=\frac{(2 p+1)P_{2 i}(u)P_{2 i+2}(u)+u\bigl(-P_{2 i+2}(u)+(u-1)P_{2 i+2 k+1}(u)\bigr)}{2 u}-\frac{P_{2 i}(u)-1}{2}\\
&=\frac{(2 p+1)P_{2 i}(u)P_{2 i+2}(u)}{2 u}-\frac{P_{2 i}(u)+P_{2 i+2}(u)-(u-1)P_{2 i+2 k+1}(u)-1}{2}\,.
\end{align*}

We also have 
\begin{align*}
&\sum_{w\in{\rm Ap}_p(A)}w^2=\sum_{y=0}^{P_{2 i}(u)/u-1}\sum_{z=0}^{u-1}\left(\left(\frac{p P_{2 i}(u)}{u}+y\right)P_{2 i+2}(u)+z P_{2 i+2 k+1}(u)\right)^2\\
&=\frac{P_{2 i}(u)}{6 u^2}\bigl(6 p^2 P_{2 i}(u)^2 P_{2 i+2}(u)^2\\
&\quad +6 p P_{2 i}(u)P_{2 i+2}(u)\bigl(P_{2 i}(u)P_{2 i+2}(u)+u\bigl(-P_{2 i+2}(u)+(u-1)P_{2 i+2 k+1}(u)\bigr)\\ 
&\quad +P_{2 i}(u)P_{2 i+2}(u)\bigl(2 P_{2 i}(u)P_{2 i+2}(u)+3 u\bigl(-P_{2 i+2}(u)+(u-1)P_{2 i+2 k+1}(u)\bigr)\\
&\quad +u^2\bigl(P_{2 i+2}(u)^2-3(u-1)P_{2 i+2}(u)P_{2 i+2 k+1}(u)+(2 u-1)(u-1)P_{2 i+2 k+1}(u)^2\bigr)\bigr)\,.  
\end{align*} 
Hence, by Lemma \ref{lem-mp} (\ref{mp-s}), we obtain 
\begin{align*}
&s_p\bigl(P_{2 i}(u),P_{2 i+2}(u),P_{2i+2 k+1}(u)\bigr)\\
&=\frac{p^2 P_{2 i}(u)^2 P_{2 i+2}(u)^2}{2 u^2}\\
&\quad +\frac{p P_{2 i}(u)P_{2 i+2}(u)}{2 u^2}\bigl(P_{2 i}(u)P_{2 i+2}(u)-u(P_{2 i}(u)+P_{2 i+2}(u)-(u-1)P_{2 i+2 k+1}(u))\bigr)\\
&\quad +\frac{1}{12 u^2}\biggl(P_{2 i}(u)P_{2 i+2}(u)\bigl(2 P_{2 i}(u)P_{2 i+2}(u)\\
&\qquad\quad -3 u\bigl(P_{2 i}(u)+P_{2 i+2}(u)-(u-1)P_{2 i+2 k+1}(u)-u)\bigr)\\
&\qquad +u^2\bigl(P_{2 i}(u)^2+P_{2 i+2}(u)^2-1\\
&\qquad\quad -(u-1)P_{2 i+2 k+1}(u)\bigl(3 P_{2 i}(u)+3 P_{2 i+2}(u)-(2 u-1)P_{2 i+2 k+1}(u))\bigr)\biggr)\,. 
\end{align*}

\begin{theorem} 
For $2\le i\le k$ and $0\le p<u(P_{2 k+1}(u)+1)/P_{2 i}(u)-1$, we have 
\begin{align*}
&n_p\bigl(P_{2 i}(u),P_{2 i+2}(u),P_{2i+2 k+1}(u)\bigr)\\
&=\frac{(2 p+1)P_{2 i}(u)P_{2 i+2}(u)}{2 u}-\frac{P_{2 i}(u)+P_{2 i+2}(u)-(u-1)P_{2 i+2 k+1}(u)-1}{2}\,,\\
&s_p\bigl(P_{2 i}(u),P_{2 i+2}(u),P_{2i+2 k+1}(u)\bigr)\\
&=\frac{p^2 P_{2 i}(u)^2 P_{2 i+2}(u)^2}{2 u^2}\\
&\quad +\frac{p P_{2 i}(u)P_{2 i+2}(u)}{2 u^2}\bigl(P_{2 i}(u)P_{2 i+2}(u)-u(P_{2 i}(u)+P_{2 i+2}(u)-(u-1)P_{2 i+2 k+1}(u))\bigr)\\
&\quad +\frac{1}{12 u^2}\biggl(P_{2 i}(u)P_{2 i+2}(u)\bigl(2 P_{2 i}(u)P_{2 i+2}(u)\\
&\qquad\quad -3 u\bigl(P_{2 i}(u)+P_{2 i+2}(u)-(u-1)P_{2 i+2 k+1}(u)-u)\bigr)\\
&\qquad +u^2\bigl(P_{2 i}(u)^2+P_{2 i+2}(u)^2-1\\
&\qquad\quad -(u-1)P_{2 i+2 k+1}(u)\bigl(3 P_{2 i}(u)+3 P_{2 i+2}(u)-(2 u-1)P_{2 i+2 k+1}(u))\bigr)\biggr)\,.  
\end{align*}
\label{th2}
\end{theorem}

\subsection{Examples}  

When $i=k=2$ and $u=2$ in Theorems \ref{th1} and \ref{th2}, we have for $p=0,1,2,3$ 
\begin{align*}
\{g_p(P_{4},P_{6},P_{9})\}_{p=0}^3&=\{1323, 1743, 2163, 2583\}\,,\\ 
\{n_p(P_{4},P_{6},P_{9})\}_{p=0}^3&=\{662, 1082, 1502, 1922\}\,,\\
\{s_p(P_{4},P_{6},P_{9})\}_{p=0}^3&=\{347209, 713239, 1255669, 1974499\}\,. 
\end{align*}
However, since $2(P_{5}+1)/P_{4}-1=4$, the formulas do not hold for $p\ge 4$.    

When $i=2$ and $k=3$, by $2(P_{7}+1)/P_{4}-1\approx 27.3$, the formulas hold for $0\le p\le 27$.

\section{The case $\{P_{2 i+1}(u),P_{2 i+3}(u),P_{2 i+k+1}(u)\}$} 

Consider the case where $A:=\{P_{2 i+1}(u),P_{2 i+3}(u),P_{2 i+k+1}(u)\}$, where $i\ge 1$ and $k\ge 3$. In this case, the structure becomes very complicated even for $p=0$. 

First of all, since $P_{2 i+1}(u)|P_{4 i+2}(u)$ and 
\begin{align*}
g\bigl(P_{2 i+1}(u),P_{2 i+3}(u)\bigr)&=P_{2 i+1}(u)P_{2 i+3}(u)-P_{2 i+1}(u)-P_{2 i+3}(u)\\
&<P_{2 i+1}(u)P_{2 i+3}(u)+P_{2 i}(u)P_{2 i+2}(u)=P_{4 i+3}(u)\,,
\end{align*}     
we have 
\begin{equation}
g\bigl(P_{2 i+1}(u),P_{2 i+3}(u),P_{2 i+k+1}(u)\bigr)=g\bigl(P_{2 i+1}(u),P_{2 i+3}(u)\bigr)\quad(k\ge 2 i+1)\,.
\label{eq:g=g}
\end{equation} 

Hereafter, we assume that $k<2 i+1$. 
In order to determine the geometrical border of the elements of ${\rm Ap}_0(A)$, consider the condition 
\begin{equation}
y P_{2 i+3}(u)\equiv z P_{2 i+k+1}(u)\pmod{P_{2 i+1}(u)}\,, 
\label{cond:1}
\end{equation} 
where 
\begin{equation}
y P_{2 i+3}(u)\not\in{\rm Ap}_0(A)\quad\hbox{and}\quad z P_{2 i+k+1}(u)\in{\rm Ap}_0(A)\,.
\label{cond:apery}
\end{equation} 
Since $P_{2 i+3}(u)\equiv u P_{2 i}(u)$ and $P_{2 i+k+1}(u)\equiv P_{2 i}(u)P_k(u)\pmod{P_{2 i+1}(u)}$, the condition (\ref{cond:1}) can be written as $u y P_{2 i}(u)\equiv z P_{2 i}(u)P_k(u)\pmod{P_{2 i+1}(u)}$. 

As $\gcd\bigl(P_{2 i}(u),P_{2 i+1}(u)\bigr)=1$, we get $u y\equiv z P_k(u)\pmod{P_{2 i+1}(u)}$. When $k<2 i+1$, we have $u y=z P_k(u)$. Therefore, if $k$ is odd, then by $P_k(u)\equiv 1\pmod{u}$, $z\equiv 0\pmod u$. The least possible positive value is $z=u$, so $y=P_k(u)$.  If $k$ is even, then by $P_k(u)\equiv 0\pmod{u}$, the least possible value is $z=1$, so $y=P_k(u)/u$. Indeed, by (\ref{pell-equ}), we get $P_k(u)P_{2 i+3}(u)-u P_{2 i+k+1}(u)=P_{2 i+1}(u)P_k(u)-u P_{2 i+1}(u)P_{k-1}(u)=P_{2 i+1}(u)P_{k-2}(u)>0$. Hence, we have $\bigl(P_k(u)/u\bigr)P_{2 i+3}(u)-P_{2 i+k+1}(u)>0$, satisfying the condition (\ref{cond:1}) and (\ref{cond:apery}). However, for example, choosing $z=y$ and $P_k(u)=u$, we have $y P_{2 i+3}(u)-z P_{2 i+k+1}(u)<0$, which against the condition (\ref{cond:apery}).  

This discussion is supported by the fact that $y P_{2 i+3}(u)\not\in{\rm Ap}_0(A)$. Namely, $y P_{2 i+3}(u)$ must be expressed in at least two different ways.  
When $z=u$ and $y=P_k(u)$ ($k$ is odd),  
$$
P_k(u)P_{2i+3}(u)=P_{k-2}(u)P_{2 i+1}(u)+u P_{2 i+k+1}(u)\,. 
$$ 
When $z=1$ and $y=P_k(u)/u$ ($k$ is even), we have 
$$
\frac{P_k(u)}{u}P_{2i+3}(u)=\frac{P_{k-2}(u)}{u}P_{2 i+1}(u)+P_{2 i+k+1}(u)\,. 
$$ 
In addition, $(P_k(u)-1)P_{2i+3}(u)$ or $(P_k(u)/u-1)P_{2i+3}(u)$ does not have another representation, respectively.

\subsection{The case $p=0$}  

First, let $k$ be odd with $k\ge 3$. 
Let $\mathfrak q$ be the largest integer with 
$$
\mathfrak q\le\frac{P_{2 i+1}(u)}{P_k(u)}\quad\hbox{and}\quad\mathfrak q\equiv 1\pmod u\,,
$$ 
and $\mathfrak r=\bigl(P_{2 i+1}(u)-\mathfrak q P_k(u)\bigr)/u$, satisfying
\begin{equation}
P_{2 i+1}(u)=\mathfrak q P_k(u)+u\mathfrak r\,. 
\label{eq:rec-k-odd}
\end{equation} 
Because of the assumption $k<2 i+1$, by $P_{2 i+1}(u)/P_k(u)>1$ such $\mathfrak q$ and $\mathfrak r$ can be chosen uniquely. In particular, $\mathfrak q$ is determined as 
$$
\mathfrak q=u\cl{\frac{1}{u}\fl{\frac{P_{2 i+1}(u)}{P_k(u)}}}-u+1\,. 
$$  
Since $\mathfrak q\ge 1$, the least possible value of $\mathfrak q$ is $\mathfrak q=1$.  
Since $P_{2 i+1}(u),\mathfrak q,P_k(u)\equiv 1\pmod u$, $\mathfrak r$ is a non-negative integer.  
From the above discussion, all the elements in ${\rm Ap}_0(A)$, where $A:=\{P_{2 i+1}(u),P_{2 i+3}(u),P_{2 i+k+1}(u)\}$, are arranged as in Table \ref{tb:ap00odd}.     

\begin{table}[htbp]
  \centering
%\scalebox{0.7}{
\begin{tabular}{ccccc}
\cline{1-2}\cline{3-4}\cline{5-5}
\multicolumn{1}{|c}{$t_{0,0}$}&$t_{1,0}$&$\cdots$&$\cdots$&\multicolumn{1}{c|}{$t_{P_{k}(u)-1,0}$}\\
\multicolumn{1}{|c}{$t_{0,1}$}&$t_{1,1}$&$\cdots$&$\cdots$&\multicolumn{1}{c|}{$t_{P_{k}(u)-1,1}$}\\
\multicolumn{1}{|c}{$\vdots$}&$\vdots$&&&\multicolumn{1}{c|}{$\vdots$}\\
\multicolumn{1}{|c}{$t_{0,\mathfrak q-1}$}&$t_{1,\mathfrak q-1}$&$\cdots$&$\cdots$&\multicolumn{1}{c|}{$t_{P_{k}(u)-1,\mathfrak q-1}$}\\
\cline{4-5}
\multicolumn{1}{|c}{$t_{0,\mathfrak q}$}&$\cdots$&\multicolumn{1}{c|}{$t_{\mathfrak r-1,\mathfrak q}$}&&\\
\multicolumn{1}{|c}{$\vdots$}&&\multicolumn{1}{c|}{$\vdots$}&&\\
\multicolumn{1}{|c}{$t_{0,\mathfrak q+u-1}$}&$\cdots$&\multicolumn{1}{c|}{$t_{\mathfrak r-1,\mathfrak q+u-1}$}&&\\
\cline{1-2}\cline{3-3}
\end{tabular}
%} 
  \caption{${\rm Ap}_0\bigl(P_{2 i+1}(u),P_{2 i+3}(u),P_{2i+k+1}(u)\bigr)$ for odd $k$}
  \label{tb:ap00odd}
\end{table}

From (\ref{eq:rec-k-odd}), the number of the elements in Table \ref{tb:ap00odd} is exactly equal to $P_{2 i+1}(u)$. 
When $0\le j\le\mathfrak q-1$, each value $t_{l,j}$ ($0\le l\le P_k(u)-1$) has only one representation but $t_{l,j}$ ($l\ge P_k(u)$) has two (or more) representations with the relation 
$$
t_{l,j+u}\equiv t_{P_k(u)+l,j}\pmod{P_{2 i+1}(u)}\quad(0\le j\le \mathfrak q-1)\,. 
$$    
When $j=\mathfrak q,\dots,\mathfrak q+u-1$, each value $t_{l,j}$ ($0\le l\le\mathfrak r-1$) has only one representation but $t_{l,j}$ ($l\ge\mathfrak r$) has two (or more) representations with the relation 
$$
t_{l,j}\equiv t_{l+\mathfrak r,j+\mathfrak q}\pmod{P_{2 i+1}(u)}\quad(0\le l\le\mathfrak r-1;\, j=0,1,\dots,u-1)\,. 
$$  
We shall show that the sequence at 
\begin{align}  
&(0,0),\dots,\bigl(P_k(u)-1,0\bigr),(0,u),\dots,\bigl(P_k(u)-1,u\bigr),\notag\\
&(0,2 u),\dots,\bigl(P_k(u)-1,2 u\bigr),\dots,(0,\mathfrak q-1),\dots,\bigl(P_k(u)-1,\mathfrak q-1\bigr),\notag\\
&(0,\mathfrak q+u-1),\dots,(\mathfrak r-1,\mathfrak q+u-1),\notag\\
&(0,u-1),\dots,\bigl(P_k(u)-1,u-1\bigr),(0,2 u-1),\dots,\bigl(P_k(u)-1,2 u-1\bigr),\dots,\notag\\
&(0,\mathfrak q-2),\dots,\bigl(P_k(u)-1,\mathfrak q-2\bigr),(0,\mathfrak q+u-2),\dots,(\mathfrak r-1,\mathfrak q+u-2),\notag\\ 
&(0,u-2),\dots,\bigl(P_k(u)-1,u-2\bigr),(0,2 u-2),\dots,\bigl(P_k(u)-1,2 u-2\bigr),\dots,\notag\\
&(0,\mathfrak q-3),\dots,\bigl(P_k(u)-1,\mathfrak q-3\bigr),(0,\mathfrak q+u-3),\dots,(\mathfrak r-1,\mathfrak q+u-3),\cdots,\notag\\ 
&(0,1),\dots,\bigl(P_k(u)-1,1\bigr),(0,u+1),\dots,\bigl(P_k(u)-1,u+1\bigr),\notag\\
&(0,2 u+1),\dots,\bigl(P_k(u)-1,2 u+1\bigr),\dots,(0,\mathfrak q-u),\dots,\bigl(P_k(u)-1,\mathfrak q-u\bigr),\notag\\
&(0,\mathfrak q),\dots,(\mathfrak r-1,\mathfrak q) 
\label{eq:seq-odd}
\end{align} 
constitutes the sequence $\{u\ell P_{2 i}(u)\}_{\ell=0}^{P_{2 i+1}(u)-1}\pmod{P_{2 i+1}(u)}$. Then by \\
$\gcd\bigl(P_{2 i}(u),P_{2 i+1}(u)\bigr)=1$ and $P_{2 i+1}(u)\equiv 1\pmod u$, the set of this sequence is equivalent to the set of the sequence $\{\ell\}_{\ell=0}^{P_{2 i+1}(u)-1}\pmod{P_{2 i+1}(u)}$. Hence, all the elements in the area of Table \ref{tb:ap00odd} constitutes a complete residue system modulo $P_{2 i+1}(u)$.        

Now, we examine the sequence (\ref{eq:seq-odd}).  
Since 
\begin{equation} 
P_{2 i+3}(u)=(u^2+1)P_{2 i+1}(u)+u P_{2 i}(u)\equiv u P_{2 i}(u)\pmod{P_{2 i+1}(u)}\,,
\label{eq:2i+3}
\end{equation} 
each row is a part of the sequence $\{u\ell P_{2 i}(u)\}_{\ell}$. Since 
\begin{equation} 
P_{2 i+k+1}(u)=P_k(u)P_{2 i}(u)+P_{k+1}(u)P_{2 i+1}(u)\,,
\label{eq:2i+k+1} 
\end{equation}  
we have for $z=0,1,\dots,\mathfrak q-1$, 
$t_{P_k(u)-1,z}\equiv\bigl((z+u)P_k(u)-u\bigr)P_{2 i}(u)$  and $t_{0,z+u}\equiv (z+u)P_k(u)P_{2 i}(u)\pmod{P_{2 i+1}(u)}$. 
In addition, by the relation (\ref{eq:rec-k-odd}), we have 
\begin{align*}
t_{\mathfrak r-1,\mathfrak q+z}&\equiv(\mathfrak r-1)\cdot u P_{2 i}(u)+(\mathfrak q+z)P_k(u)P_{2 i}(u)\\
&\equiv\bigl(z P_k(u)-u\bigr)P_{2 i}(u)\pmod{P_{2 i+1}(u)}
\end{align*} 
and $t_{0,z}\equiv z P_k(u)P_{2 i}(u)\pmod{P_{2 i+1}(u)}$.  
Finally, the final element is 
$$
t_{\mathfrak r-1,\mathfrak q}\equiv(\mathfrak r-1)\cdot u P_{2 i}(u)+\mathfrak q P_k(u)P_{2 i}(u)\equiv -u P_{2 i}(u)\pmod{P_{2 i+1}(u)}\,. 
$$ 
From Table \ref{tb:ap00odd}, we see that the largest element in ${\rm Ap}_0(A)$ is $t_{\mathfrak r-1,\mathfrak q+u-1}$ or $t_{P_k(u)-1,\mathfrak q-1}$.  
The condition $t_{\mathfrak r-1,\mathfrak q+u-1}>t_{P_k(u)-1,\mathfrak q-1}$ is equivalent to $u\mathfrak r P_{2 i}(u)>(P_{k-2}(u)-(u^2+1)\mathfrak r)P_{2 i+1}(u)$. 

Therefore, if $u\mathfrak r P_{2 i}(u)>(P_{k-2}(u)-(u^2+1)\mathfrak r)P_{2 i+1}(u)$, then 
$$
g_0(A)=(\mathfrak r-1)P_{2 i+3}(u)+(\mathfrak q+u-1)P_{2 i+k+1}(u)-P_{2 i+1}(u)\,.  
$$   
Otherwise,  
$$
g_0(A)=\bigl(P_k(u)-1\bigr)P_{2 i+3}(u)+(\mathfrak q-1)P_{2 i+k+1}(u)-P_{2 i+1}(u)\,.  
$$ 
\bigskip

Next, let $k$ be even with $k\ge 4$. All the elements in ${\rm Ap}_0(A)$, where $A:=\{P_{2 i+1}(u),P_{2 i+3}(u),P_{2 i+k+1}(u)\}$, are arranged as in Table \ref{tb:ap00even}. 
Let $q$ be the largest integer not exceeding $u P_{2 i+1}/P_k$, and $r=P_{2 i+1}-q\bigl(P_k(u)/u\bigr)$, satisfying
\begin{equation}
P_{2 i+1}(u)=q\frac{P_k(u)}{u}+r\,. 
\label{eq:rec-k-even}
\end{equation} 
From the above discussion, all the elements in ${\rm Ap}_0(A)$, where \\ 
$A:=\{P_{2 i+1}(u),P_{2 i+3}(u),P_{2 i+k+1}(u)\}$, are arranged as in Table \ref{tb:ap00even}.

\begin{table}[htbp]
  \centering
%\scalebox{0.7}{
\begin{tabular}{ccccc}
\cline{1-2}\cline{3-4}\cline{5-5}
\multicolumn{1}{|c}{$t_{0,0}$}&$t_{1,0}$&$\cdots$&$\cdots$&\multicolumn{1}{c|}{$t_{P_{k}(u)/u-1,0}$}\\
\multicolumn{1}{|c}{$t_{0,1}$}&$t_{1,1}$&$\cdots$&$\cdots$&\multicolumn{1}{c|}{$t_{P_{k}(u)/u-1,1}$}\\
\multicolumn{1}{|c}{$\vdots$}&$\vdots$&&&\multicolumn{1}{c|}{$\vdots$}\\
\multicolumn{1}{|c}{$t_{0,q-1}$}&$t_{1,q-1}$&$\cdots$&$\cdots$&\multicolumn{1}{c|}{$t_{P_{k}(u)/u-1,q-1}$}\\
\cline{4-5}
\multicolumn{1}{|c}{$t_{0,q}$}&$\cdots$&\multicolumn{1}{c|}{$t_{r-1,q}$}&&\\
\cline{1-2}\cline{3-3}
\end{tabular}
%} 
  \caption{${\rm Ap}_0\bigl(P_{2 i+1}(u),P_{2 i+3}(u),P_{2i+k+1}(u)\bigr)$ for even $k$}
  \label{tb:ap00even}
\end{table}

Similarly to the case where $k$ is odd, the sequence 
\begin{align*}
&(0,0),\dots,\bigl(P_k(u)/u-1,0\bigr),(0,1),\dots,\bigl(P_k(u)/u-1,1\bigr),\\
&(0,2),\dots,\bigl(P_k(u)/u-1,2\bigr),\dots,(0,q-1),\dots,\bigl(P_k(u)/u-1,q-1\bigr),\\
&(0,q),\dots,(r-1,q)
\end{align*}
constitutes the sequence $\{u\ell P_{2 i}(u)\}_{\ell=0}^{P_{2 i+1}(u)-1}\pmod{P_{2 i+1}(u)}$, so $\{\ell\}_{\ell=0}^{P_{2 i+1}(u)-1}\pmod{P_{2 i+1}(u)}$. Hence, all the elements in the area of Table \ref{tb:ap00even} constitutes a complete residue system modulo $P_{2 i+1}$.   
From Table \ref{tb:ap00even}, the largest element in ${\rm Ap}_0(A)$ is $t_{r-1,q}$ or $t_{P_k(u)/u-1,q-1}$.  
The condition $t_{r-1,q}>t_{P_k(u)/u-1,q-1}$ is equivalent to $u r P_{2 i}(u)>\bigl(P_{k-2}(u)/u-(u^2+1)r\bigr)P_{2 i+1}(u)$.  
Therefore, if $u r P_{2 i}(u)>\bigl(P_{k-2}(u)/u-(u^2+1)r\bigr)P_{2 i+1}(u)$, then 
$$
g_0(A)=(r-1)P_{2 i+3}(u)+q P_{2 i+k+1}(u)-P_{2 i+1}(u)\,.  
$$  
Otherwise, 
$$
g_0(A)=\left(\frac{P_k(u)}{u}-1\right)P_{2 i+3}(u)+(q-1)P_{2 i+k+1}(u)-P_{2 i+1}(u)\,.  
$$

\subsection{Examples}  
\label{examples4}

Consider the case where $k$ is odd. Let $u=2$.  
For $A:=\{P_9,P_{11},P_{12}\}$, by $197<P_9/P_5<198$, we get  
$$
\mathfrak q=197\quad\hbox{and}\quad \mathfrak r=\frac{P_9-197 P_5}{2}=0\,. 
$$   
Hence, by $2\cdot 0\cdot P_{6}<(P_{3}-5\cdot 0)P_{7}$, we have 
$$
g_0(A)=(P_3-1)P_{11}+(197-1)P_{12}-P_{9}=2738539\,.  
$$ 

For $A:=\{P_{11},P_{13},P_{14}\}$, by $1148<P_{11}/P_3<1149$, we get  
$$
\mathfrak q=1147\quad\hbox{and}\quad \mathfrak r=\frac{P_{11}-1147 P_3}{2}=3 
$$   
(Note that $q$ must be odd).   
Hence, by $2\cdot 3\cdot P_{10}>(P_{1}-5\cdot 3)P_{11}$, we have 
$$
g_0(A)=(3-1)P_{13}+(1147+1)P_{14}-P_{11}=92798917\,.  
$$ 
\bigskip 

Consider the case where $k$ is even.   
For $A:=\{P_{7},P_{9},P_{11}\}$, by $28<2 P_{7}/P_4<29$, we get   
$$
q=28\quad\hbox{and}\quad r=P_{7}-28\frac{P_4}{2}=1\,.  
$$    
Hence, by $2\cdot 1\cdot P_{6}>(P_{2}/2-5\cdot 1)P_{7}$, we have 
$$
g_0(A)=(1-1)P_{9}+28 P_{11}-P_{7}=160579\,.  
$$  

For $A:=\{P_{11},P_{13},P_{17}\}$, by $164<2 P_{11}/P_6<165$, we get  
$$
q=164\quad\hbox{and}\quad r=P_{11}-164\frac{P_6}{2}=1\,.  
$$    
Hence, by $2\cdot 1\cdot P_{10}<(P_{4}/2-5\cdot 1)P_{11}$, we have 
$$
g_0(A)=\left(\frac{P_6}{2}-1\right)P_{13}+(164-1)P_{17}-P_{11}=186412240\,.  
$$

\subsection{Odd $k$}  

First, let $k$ be odd for $A:=\{P_{2 i+1}(u),P_{2 i+3}(u),P_{2i+k+1}(u)\}$.  

\subsubsection{The case $p=1$}  

Geometrically speaking, the position of each element of ${\rm Ap}_1(A)$ is determined from the position of an element of ${\rm Ap}_0(A)$. Namely, the most part except for the first two rows in ${\rm Ap}_0(A)$ is shifted up by two rows as it is and placed immediately to the right of ${\rm Ap}_0(A)$. Only the first two rows of ${\rm Ap}_0(A)$ are arranged to stick together at the end of ${\rm Ap}_0(A)$, and the remaining two rows are further arranged from the left at the end. The elements of ${\rm Ap}_1(A)$ are arranged in three parts (two parts when $\mathfrak r=0$). See Table \ref{tb:ap01odd}.

\begin{table}[htbp]
  \centering
\scalebox{0.6}{
\begin{tabular}{cccccccccc}
\cline{1-2}\cline{3-4}\cline{5-6}\cline{7-8}\cline{9-10}
\multicolumn{1}{|c}{$t_{0,0}$}&$t_{1,0}$&$\cdots$&$\cdots$&\multicolumn{1}{c|}{$t_{P_{k}(u)-1,0}$}&$t_{P_{k}(u),0}$&$t_{P_{k}(u)+1,0}$&$\cdots$&$\cdots$&\multicolumn{1}{c|}{$t_{2 P_{k}(u)-1,0}$}\\
\multicolumn{1}{|c}{$t_{0,1}$}&$t_{1,1}$&$\cdots$&$\cdots$&\multicolumn{1}{c|}{$t_{P_{k}(u)-1,1}$}&$t_{P_{k}(u),1}$&$t_{P_{k}(u)+1,1}$&$\cdots$&$\cdots$&\multicolumn{1}{c|}{$t_{2 P_{k}(u)-1,1}$}\\
\multicolumn{1}{|c}{$\vdots$}&$\vdots$&&&\multicolumn{1}{c|}{$\vdots$}&$\vdots$&$\vdots$&&&\multicolumn{1}{c|}{$\vdots$}\\
\multicolumn{1}{|c}{$t_{0,\mathfrak q-u-1}$}&$t_{1,\mathfrak q-u-1}$&$\cdots$&$\cdots$&\multicolumn{1}{c|}{$t_{P_{k}(u)-1,\mathfrak q-u-1}$}&$t_{P_{k}(u),\mathfrak q-u-1}$&$t_{P_{k}(u)+1,\mathfrak q-u-1}$&$\cdots$&$\cdots$&\multicolumn{1}{c|}{$t_{2 P_{k}(u)-1,\mathfrak q-u-1}$}\\
\cline{9-10} 
\multicolumn{1}{|c}{$t_{0,\mathfrak q-u}$}&$t_{1,\mathfrak q-u}$&$\cdots$&$\cdots$&\multicolumn{1}{c|}{$t_{P_{k}(u)-1,\mathfrak q-u}$}&$t_{P_{k}(u),\mathfrak q-u}$&$\cdots$&\multicolumn{1}{c|}{$t_{P_{k}(u)+\mathfrak r-1,\mathfrak q-u}$}&&\\
\multicolumn{1}{|c}{$\vdots$}&$\vdots$&&&\multicolumn{1}{c|}{$\vdots$}&$\vdots$&&\multicolumn{1}{c|}{$\vdots$}&&\\
\multicolumn{1}{|c}{$t_{0,\mathfrak q-1}$}&$t_{1,\mathfrak q-1}$&$\cdots$&$\cdots$&\multicolumn{1}{c|}{$t_{P_{k}(u)-1,\mathfrak q-1}$}&$t_{P_{k}(u),\mathfrak q-1}$&$\cdots$&\multicolumn{1}{c|}{$t_{P_{k}(u)+\mathfrak r-1,\mathfrak q-1}$}&&\\
\cline{4-4}\cline{5-6}\cline{7-8}
\multicolumn{1}{|c}{$t_{0,\mathfrak q}$}&$\cdots$&\multicolumn{1}{c|}{$t_{\mathfrak r-1,\mathfrak q}$}&$\cdots$&\multicolumn{1}{c|}{$t_{P_{k}(u)-1,\mathfrak q}$}&&&&&\\
\multicolumn{1}{|c}{$\vdots$}&&\multicolumn{1}{c|}{$\vdots$}&&\multicolumn{1}{c|}{$\vdots$}&&&&&\\
\multicolumn{1}{|c}{$t_{0,q+\mathfrak u-1}$}&$\cdots$&\multicolumn{1}{c|}{$t_{\mathfrak r-1,\mathfrak q+u-1}$}&$\cdots$&\multicolumn{1}{c|}{$t_{P_{k}(u)-1,\mathfrak q+u-1}$}&&&&&\\
\cline{1-2}\cline{3-4}\cline{5-5}
\multicolumn{1}{|c}{$t_{0,\mathfrak q+u}$}&$\cdots$&\multicolumn{1}{c|}{$t_{\mathfrak r-1,\mathfrak q+u}$}&&&&&&&\\
\multicolumn{1}{|c}{$\vdots$}&&\multicolumn{1}{c|}{$\vdots$}&&&&&&&\\
\multicolumn{1}{|c}{$t_{0,\mathfrak q+2 u-1}$}&$\cdots$&\multicolumn{1}{c|}{$t_{\mathfrak r-1,\mathfrak q+2 u-1}$}&&&&&&&\\
\cline{1-2}\cline{3-3}
\end{tabular}
} 
  \caption{${\rm Ap}_p\bigl(P_{2 i+1}(u),P_{2 i+3}(u),P_{2i+k+1}(u)\bigr)$ ($p=0,1$) for odd $k$}
  \label{tb:ap01odd}
\end{table}

These geometrical interpretations are supported by the congruence relations:  
\begin{align*}
&t_{y,z+u}\equiv t_{P_k(u)+y,z}\pmod{P_{2 i+1}(u)}\quad(0\le y\le P_k(u)-1,\,0\le z\le\mathfrak q-1)\,,\\
&t_{y,z}\equiv t_{\mathfrak r+y,\mathfrak q+z}\pmod{P_{2 i+1}(u)}\quad(0\le y\le P_k(u)-\mathfrak r-1,\,0\le z\le u-1)\,,\\
&t_{P_k(u)+y-\mathfrak r,z}\equiv t_{y,\mathfrak q+z+u}\pmod{P_{2 i+1}(u)}\quad(0\le y\le\mathfrak r-1,\,0\le z\le u-1)\,
\end{align*}
In addition, all the elements in ${\rm Ap}_1(A)$ are represented in at least two ways because 
\begin{align*}
&(P_k(u)+y)P_{2 i+3}(u)+z P_{2 i+k+1}(u)\\
&=P_{k-2}(u)P_{2 i+1}(u)+y P_{2 i+3}(u)+(z+u)P_{2 i+k+1}(u)\,,\\
&(\mathfrak r+y)P_{2 i+3}(u)+(\mathfrak q+z)P_{2 i+k+1}(u)\\
&=\bigl(P_{2 i}(u)+\mathfrak q P_{k+1}(u)+(u^2+1)\mathfrak r\bigr)P_{2 i+1}(u)+y P_{2 i+3}(u)+z P_{2 i+k+1}(u)\,,\\ 
&y P_{2 i+3}(u)+(\mathfrak q+z+u)P_{2 i+k+1}(u)\\
&=(P_{2 i}(u)+\mathfrak q P_{k+1}(u)-P_{k-2}(u)+(u^2+1)\mathfrak r\bigr)P_{2 i+1}(u)\\
&\quad +\bigl(P_k(u)+y-\mathfrak r)P_{2 i+3}(u)+z P_{2 i+k+1}(u)\,. 
\end{align*}

There are four candidates to take the largest value of ${\rm Ap}_1(A)$: 
$$
t_{\mathfrak r-1,\mathfrak q+2 u-1},\quad t_{P_k(u)-1,\mathfrak q+u-1},\quad t_{P_k(u)+\mathfrak r-1,\mathfrak q-1},\quad t_{2 P_k(u)-1,\mathfrak q-u-1}\,. 
$$   
Since $2 u P_{2 i+k+1}(u)>P_k(u)P_{2 i+3}(u)$, we see that 
$$
t_{\mathfrak r-1,\mathfrak q+2 u-1}>t_{P_k(u)+\mathfrak r-1,\mathfrak q-1}\quad\hbox{and}\quad t_{P_k(u)-1,\mathfrak q+u-1}>t_{2 P_k(u)-1,\mathfrak q-2 u-1}\,. 
$$ 
Therefore, if $u\mathfrak r P_{2 i}(u)>(P_{k-2}(u)-(u^2+1)\mathfrak r)P_{2 i+1}(u)$, then 
$$
g_1(A)=(\mathfrak r-1)P_{2 i+3}(u)+(\mathfrak q+2 u-1)P_{2 i+k+1}(u)-P_{2 i+1}(u)\,.  
$$   
Otherwise,  
$$
g_1(A)=\bigl(P_k(u)-1\bigr)P_{2 i+3}(u)+(\mathfrak q+u-1)P_{2 i+k+1}(u)-P_{2 i+1}(u)\,.  
$$

\subsubsection{The case $p\ge 2$}  

When $p=2$, the position of the elements of ${\rm Ap}_2(A)$ are similarly determined by that of ${\rm Ap}_1(A)$. Namely, the most part except for the first two rows in ${\rm Ap}_1(A)$ is shifted up by two rows as it is and placed immediately to the right of ${\rm Ap}_1(A)$. Only the first two rows of ${\rm Ap}_1(A)$ are arranged to stick together at the end of ${\rm Ap}_0(A)$, and the remaining two rows are further arranged from the left at the end. See Table \ref{tb:ap02odd}.

\begin{table}[htbp]
  \centering
\scalebox{0.5}{
\begin{tabular}{ccccccccccccccc}
\cline{1-2}\cline{3-4}\cline{5-6}\cline{7-8}\cline{9-10}\cline{11-12}\cline{13-14}\cline{15-15}
\multicolumn{1}{|c}{}&&&&\multicolumn{1}{c|}{}&&&&&\multicolumn{1}{c|}{}&&&&&\multicolumn{1}{c|}{}\\ 
\multicolumn{1}{|c}{}&&&&\multicolumn{1}{c|}{}&&&&&\multicolumn{1}{c|}{}&&&&$\cdots$&\multicolumn{1}{c|}{$t_{3 P_k(u)-1,\mathfrak q-2 u-1}$}\\ 
\cline{14-15}
\multicolumn{1}{|c}{}&&&&\multicolumn{1}{c|}{}&&&&&\multicolumn{1}{c|}{}&&&\multicolumn{1}{c|}{}&&\\ 
\multicolumn{1}{|c}{}&&&&\multicolumn{1}{c|}{}&&&&&\multicolumn{1}{c|}{}&&$\cdots$&\multicolumn{1}{c|}{$t_{2 P_k(u)+\mathfrak r-1,\mathfrak q-u-1}$}&&\\ 
\cline{9-10}\cline{11-12}\cline{13-13}
\multicolumn{1}{|c}{}&&&&\multicolumn{1}{c|}{}&&&\multicolumn{1}{c|}{}&&\multicolumn{1}{c|}{}&&&&&\\ 
\multicolumn{1}{|c}{}&&&&\multicolumn{1}{c|}{}&&&\multicolumn{1}{c|}{}&$\cdots$&\multicolumn{1}{c|}{$t_{2 P_k(u)-1,\mathfrak q-1}$}&&&&&\\ 
\cline{4-4}\cline{5-6}\cline{7-8}\cline{9-10}
\multicolumn{1}{|c}{}&&\multicolumn{1}{c|}{}&&\multicolumn{1}{c|}{}&&&\multicolumn{1}{c|}{}&&&&&&&\\ 
\multicolumn{1}{|c}{}&&\multicolumn{1}{c|}{}&&\multicolumn{1}{c|}{}&&$\cdots$&\multicolumn{1}{c|}{$t_{P_k(u)+\mathfrak r-1,\mathfrak q+u-1}$}&&&&&&&\\ 
\cline{1-2}\cline{3-4}\cline{5-6}\cline{7-8}
\multicolumn{1}{|c}{}&&\multicolumn{1}{c|}{}&&\multicolumn{1}{c|}{}&&&&&&&&&&\\ 
\multicolumn{1}{|c}{}&&\multicolumn{1}{c|}{}&$\cdots$&\multicolumn{1}{c|}{$t_{P_k(u)-1,\mathfrak q+2 u-1}$}&&&&&&&&&&\\ 
\cline{1-2}\cline{3-4}\cline{5-5}
\multicolumn{1}{|c}{}&&\multicolumn{1}{c|}{}&&&&&&&&&&&&\\ 
\multicolumn{1}{|c}{}&$\cdots$&\multicolumn{1}{c|}{$t_{\mathfrak r-1,\mathfrak q+3 u-1}$}&&&&&&&&&&&&\\ 
\cline{1-2}\cline{3-3}
\end{tabular}
} 
  \caption{${\rm Ap}_p\bigl(P_{2 i+1}(u),P_{2 i+3}(u),P_{2i+k+1}(u)\bigr)$ ($p=0,1,2$) for odd $k$}
  \label{tb:ap02odd}
\end{table}

Geometrically, there seem to be six candidates to take the largest value of ${\rm Ap}_2(A)$, but because of 
\begin{align*}
&t_{\mathfrak r-1,\mathfrak q+3 u-1}>t_{P_k(u)+\mathfrak r-1,\mathfrak q+u-1}>t_{2 P_k(u)+\mathfrak r-1,\mathfrak q-u-1}\,,\\
&t_{P_k(u)-1,\mathfrak q+2 u-1}>t_{2 P_k(u)-1,\mathfrak q-1}>t_{3 P_k(u)-1,\mathfrak q-2 u-1}\,, 
\end{align*}
there are only two possibilities: $t_{\mathfrak r-1,\mathfrak q+3 u-1}$ or $t_{P_k(u)-1,\mathfrak q+2 u-1}$.  
Therefore, if $u\mathfrak r P_{2 i}(u)>\bigl(P_{k-2}(u)-(u^2+1)\mathfrak r\bigr)P_{2 i+1}(u)$, then 
$$
g_2(A)=(\mathfrak r-1)P_{2 i+3}(u)+(\mathfrak q+3 u-1)P_{2 i+k+1}(u)-P_{2 i+1}(u)\,.  
$$   
Otherwise,  
$$
g_2(A)=\bigl(P_k(u)-1\bigr)P_{2 i+3}(u)+(\mathfrak q+2 u-1)P_{2 i+k+1}(u)-P_{2 i+1}(u)\,.  
$$ 
\bigskip

When $p=3,4,\dots$, similar discussions can be applied. All the elements of ${\rm Ap}_p(A)$ are obtained from those of ${\rm Ap}_{p-1}(A)$.  The positions of the elements of ${\rm Ap}_p(A)$ below the left-most block and the positions of ${\rm Ap}_p(A)$ in the right-most block are arranged as shown in Table \ref{tb:ap0podd}.  

\begin{table}[htbp]
  \centering
\scalebox{0.5}{
\begin{tabular}{ccccccccccccccc}
\cline{11-12}\cline{13-14}\cline{15-15}
&&&&&&&&&\multicolumn{1}{c|}{}&&&&&\multicolumn{1}{c|}{}\\ 
&&&&&&&&&\multicolumn{1}{c|}{}&&&&$\cdots$&\multicolumn{1}{c|}{$t_{(p+1)P_k(u)-1,\mathfrak q-p u-1}$}\\ 
\cline{14-15}
&&&&&&&&&\multicolumn{1}{c|}{}&&&\multicolumn{1}{c|}{}&&\\ 
&&&&&&&&$\cdots$&\multicolumn{1}{c|}{$\cdots$}&&$\cdots$&\multicolumn{1}{c|}{$t_{p P_k(u)+\mathfrak r-1,\mathfrak q-(p-1)u-1}$}&&\\ 
\cline{11-12}\cline{13-13}
&&&&&&&&&&&&&&\\ 
&&&&&$\iddots$&&&&&&&&&\\ 
\cline{4-5}
&&\multicolumn{1}{c|}{}&&\multicolumn{1}{c|}{}&&&&&&&&&&\\ 
$\vdots$&&\multicolumn{1}{c|}{}&$\cdots$&\multicolumn{1}{c|}{$t_{P_k(u)-1,\mathfrak q+p u-1}$}&&&&&&&&&&\\ 
\cline{1-2}\cline{3-4}\cline{5-5}
\multicolumn{1}{|c}{}&&\multicolumn{1}{c|}{}&&&&&&&&&&&&\\ 
\multicolumn{1}{|c}{}&$\cdots$&\multicolumn{1}{c|}{$t_{\mathfrak r-1,\mathfrak q+(p+1)u-1}$}&&&&&&&&&&&&\\ 
\cline{1-2}\cline{3-3}
\end{tabular}
} 
  \caption{${\rm Ap}_p\bigl(P_{2 i+1}(u),P_{2 i+3}(u),P_{2i+k+1}(u)\bigr)$ for odd $k$}
  \label{tb:ap0podd}
\end{table}

This situation is continued as long as $z=\mathfrak q-p u-1\ge 0$. However, when $p=(\mathfrak q+u-1)/u$, the shape of the block on the right side collapses. Thus, the regularity of taking the maximum value of ${\rm Ap}_p(A)$ is broken. 
%However, when $p=(\mathfrak q+3)/2$, the element coming to the bottom left also collapses, and the regularity of taking the maximum value is broken.

\begin{theorem}  
Let $i\ge 1$ and $k$ be odd with $3\le k\le 2 i-1$.  
Let $\mathfrak q$ be the largest integer not exceeding $P_{2 i+1}(u)/P_k(u)$ with $\mathfrak q\equiv 1\pmod{u}$, and $\mathfrak r=\bigl(P_{2 i+1}(u)-\mathfrak q P_k(u)\bigr)/2$. Then for $0\le p\le(\mathfrak q-1)/u$, if $u\mathfrak r P_{2 i}(u)>\bigl(P_{k-2}(u)-(u^2+1)\mathfrak r\bigr)P_{2 i+1}(u)$, then 
\begin{multline*}
g_p\bigl(P_{2 i+1}(u),P_{2 i+3}(u),P_{2 i+k+1}(u)\bigr)\\
=(\mathfrak r-1)P_{2 i+3}(u)+\bigl(\mathfrak q+(p+1)u-1\bigr)P_{2 i+k+1}(u)-P_{2 i+1}(u)\,.  
\end{multline*}  
Otherwise,  
\begin{multline*}
g_p\bigl(P_{2 i+1}(u),P_{2 i+3}(u),P_{2 i+k+1}(u)\bigr)\\
=\bigl(P_k(u)-1\bigr)P_{2 i+3}(u)+(\mathfrak q+p u-1)P_{2 i+k+1}(u)-P_{2 i+1}(u)\,.  
\end{multline*}
\label{th:i-k-odd}
\end{theorem}

\subsection{Even $k$}

Next, let $k$ be even for $A:=\{P_{2 i+1}(u),P_{2 i+3}(u),P_{2i+k+1}(u)\}$.  

Geometrically speaking, the position of each element of ${\rm Ap}_1(A)$ is determined from the position of an element of ${\rm Ap}_0(A)$. Namely, the most part except for the first row in ${\rm Ap}_0(A)$ is shifted up by one row as it is and placed immediately to the right of ${\rm Ap}_0(A)$. Only the first row of ${\rm Ap}_0(A)$ are arranged to stick together at the end of ${\rm Ap}_0(A)$, and the remaining one row is further arranged from the left at the end. The elements of ${\rm Ap}_1(A)$ are arranged in three parts (two parts when $r=0$). See Table \ref{tb:ap01even}.

\begin{table}[htbp]
  \centering
\scalebox{0.6}{
\begin{tabular}{cccccccccc}
\cline{1-2}\cline{3-4}\cline{5-6}\cline{7-8}\cline{9-10}
\multicolumn{1}{|c}{$t_{0,0}$}&$t_{1,0}$&$\cdots$&$\cdots$&\multicolumn{1}{c|}{$t_{P_{k}(u)/u-1,0}$}&$t_{P_{k}(u)/u,0}$&$t_{P_{k}(u)/u+1,0}$&$\cdots$&$\cdots$&\multicolumn{1}{c|}{$t_{2 P_{k}(u)/u-1,0}$}\\
\multicolumn{1}{|c}{$t_{0,1}$}&$t_{1,1}$&$\cdots$&$\cdots$&\multicolumn{1}{c|}{$t_{P_{k}(u)/u-1,1}$}&$t_{P_{k}(u)/u,1}$&$t_{P_{k}(u)/u+1,1}$&$\cdots$&$\cdots$&\multicolumn{1}{c|}{$t_{2 P_{k}(u)/u-1,1}$}\\
\multicolumn{1}{|c}{$\vdots$}&$\vdots$&&&\multicolumn{1}{c|}{$\vdots$}&$\vdots$&$\vdots$&&&\multicolumn{1}{c|}{$\vdots$}\\
\multicolumn{1}{|c}{$t_{0,q-2}$}&$t_{1,q-2}$&$\cdots$&$\cdots$&\multicolumn{1}{c|}{$t_{P_{k}(u)/u-1,q-2}$}&$t_{P_{k}(u)/u,q-2}$&$t_{P_{k}(u)/u+1,q-2}$&$\cdots$&$\cdots$&\multicolumn{1}{c|}{$t_{2 P_{k}(u)/u-1,q-2}$}\\
\cline{9-10} 
\multicolumn{1}{|c}{$t_{0,q-1}$}&$t_{1,q-1}$&$\cdots$&$\cdots$&\multicolumn{1}{c|}{$t_{P_{k}(u)/u-1,q-1}$}&$t_{P_{k}(u)/u,\mathfrak q-1}$&$\cdots$&\multicolumn{1}{c|}{$t_{P_{k}(u)/u+r-1,q-1}$}&&\\
\cline{4-4}\cline{5-6}\cline{7-8}
\multicolumn{1}{|c}{$t_{0,q}$}&$\cdots$&\multicolumn{1}{c|}{$t_{r-1,q}$}&$\cdots$&\multicolumn{1}{c|}{$t_{P_{k}(u)/u-1,q}$}&&&&&\\
\cline{1-2}\cline{3-4}\cline{5-5}
\multicolumn{1}{|c}{$t_{0,q+1}$}&$\cdots$&\multicolumn{1}{c|}{$t_{r-1,q+1}$}&&&&&&&\\
\cline{1-2}\cline{3-3}
\end{tabular}
} 
  \caption{${\rm Ap}_p\bigl(P_{2 i+1}(u),P_{2 i+3}(u),P_{2i+k+1}(u)\bigr)$ ($p=0,1$) for even $k$}
  \label{tb:ap01even}
\end{table}

When $p=2$, the position of the elements of ${\rm Ap}_2(A)$ are similarly determined by that of ${\rm Ap}_1(A)$. Namely, the most part except for the first two rows in ${\rm Ap}_1(A)$ is shifted up by one row as it is and placed immediately to the right of ${\rm Ap}_1(A)$. Only the first row of ${\rm Ap}_1(A)$ are arranged to stick together at the end of ${\rm Ap}_0(A)$, and the remaining two rows are further arranged from the left at the end. See Table \ref{tb:ap02even}.

\begin{table}[htbp]
  \centering
\scalebox{0.5}{
\begin{tabular}{ccccccccccccccc}
\cline{1-2}\cline{3-4}\cline{5-6}\cline{7-8}\cline{9-10}\cline{11-12}\cline{13-14}\cline{15-15}
\multicolumn{1}{|c}{}&&&&\multicolumn{1}{c|}{}&&&&&\multicolumn{1}{c|}{}&&&&&\multicolumn{1}{c|}{}\\ 
\multicolumn{1}{|c}{}&&&&\multicolumn{1}{c|}{}&&&&&\multicolumn{1}{c|}{}&&&&$\cdots$&\multicolumn{1}{c|}{$t_{3 P_k(u)/u-1,q-3}$}\\ 
\cline{14-15} 
\multicolumn{1}{|c}{}&&&&\multicolumn{1}{c|}{}&&&&&\multicolumn{1}{c|}{}&&$\cdots$&\multicolumn{1}{c|}{$t_{2 P_k(u)/u+r-1,q-2}$}&&\\ 
\cline{9-10}\cline{11-12}\cline{13-13}
\multicolumn{1}{|c}{}&&&&\multicolumn{1}{c|}{}&&&\multicolumn{1}{c|}{}&$\cdots$&\multicolumn{1}{c|}{$t_{2 P_k(u)/u-1,q-1}$}&&&&&\\ 
\cline{4-4}\cline{5-6}\cline{7-8}\cline{9-10}
\multicolumn{1}{|c}{}&&\multicolumn{1}{c|}{}&&\multicolumn{1}{c|}{}&&$\cdots$&\multicolumn{1}{c|}{$t_{P_k(u)/u+r-1,q}$}&&&&&&&\\ 
\cline{1-2}\cline{3-4}\cline{5-6}\cline{7-8} 
\multicolumn{1}{|c}{}&&\multicolumn{1}{c|}{}&$\cdots$&\multicolumn{1}{c|}{$t_{P_k(u)/u-1,q+1}$}&&&&&&&&&&\\ 
\cline{1-2}\cline{3-4}\cline{5-5}
\multicolumn{1}{|c}{}&$\cdots$&\multicolumn{1}{c|}{$t_{\mathfrak r-1,q+2}$}&&&&&&&&&&&&\\ 
\cline{1-2}\cline{3-3}
\end{tabular}
} 
  \caption{${\rm Ap}_p\bigl(P_{2 i+1}(u),P_{2 i+3}(u),P_{2i+k+1}(u)\bigr)$ ($p=0,1,2$) for even $k$}
  \label{tb:ap02even}
\end{table}

When $p=3,4,\dots$, similar discussions can be applied. All the elements of ${\rm Ap}_p(A)$ are obtained from those of ${\rm Ap}_{p-1}(A)$.  The positions of the elements of ${\rm Ap}_p(A)$ below the left-most block and the positions of ${\rm Ap}_p(A)$ in the right-most block are arranged as shown in Table \ref{tb:ap0peven}.  

\begin{table}[htbp]
  \centering
\scalebox{0.6}{
\begin{tabular}{ccccccccccccccc}
\cline{11-12}\cline{13-14}\cline{15-15}
&&&&&&&&&\multicolumn{1}{c|}{}&&&&&\multicolumn{1}{c|}{}\\ 
&&&&&&&&&\multicolumn{1}{c|}{}&&&&$\cdots$&\multicolumn{1}{c|}{$t_{(p+1)P_k(u)/u-1,q-p-1}$}\\ 
\cline{14-15}
&&&&&&&&$\cdots$&\multicolumn{1}{c|}{$\cdots$}&&$\cdots$&\multicolumn{1}{c|}{$t_{p P_k(u)/u+r-1,q-p}$}&&\\ 
\cline{11-12}\cline{13-13}
&&&&&$\iddots$&&&&&&&&&\\ 
\cline{4-5} 
$\vdots$&&\multicolumn{1}{c|}{}&$\cdots$&\multicolumn{1}{c|}{$t_{P_k(u)/u-1,q+p-1}$}&&&&&&&&&&\\ 
\cline{1-2}\cline{3-4}\cline{5-5}
\multicolumn{1}{|c}{}&$\cdots$&\multicolumn{1}{c|}{$t_{\mathfrak r-1,q+p}$}&&&&&&&&&&&&\\ 
\cline{1-2}\cline{3-3}
\end{tabular}
} 
  \caption{${\rm Ap}_p\bigl(P_{2 i+1}(u),P_{2 i+3}(u),P_{2i+k+1}(u)\bigr)$ for even $k$}
  \label{tb:ap0peven}
\end{table} 

As we can see that the largest element in ${\rm Ap}_p(A)$ is either $t_{r-1,q+p}$ or $t_{P_k(u)/u-1,q+p-1}$, if $u r P_{2 i}(u)>\bigl(P_{k-2}(u)/u-(u^2+1)r\bigr)P_{2 i+1}(u)$, then 
$$
g_p(A)=(r-1)P_{2 i+3}(u)+(q+p)P_{2 i+k+1}(u)-P_{2 i+1}(u)\,. 
$$  
Otherwise, 
$$
g_p(A)=\left(\frac{P_k(u)}{u}-1\right)P_{2 i+3}(u)+(q+p-1)P_{2 i+k+1}(u)-P_{2 i+1}(u)\,. 
$$  
This situation is regularly continued until $q-p-1\ge 0$, as seen in the right-most block in Table \ref{tb:ap0peven}. When $p=q$, only the rightmost block is shortened to the side, but the other parts appear in accordance with the rules because there is only one (or zero) line as a remainder block. Therefore, the regularity of taking the maximum value of ${\rm Ap}_p(A)$ is still maintained. However, when $p=q+1$, the element coming to the bottom left becomes shorter, and the regularity of taking the maximum value is broken.  

\begin{theorem} 
Let $i\ge 1$ and $k$ be even with $4\le k\le 2 i$.   
Let $q=\fl{u P_{2 i+1}(u)/P_k(u)}$ and $r=P_{2 i+1}(u)-q\bigl(P_k(u)/u\bigr)$. Then for $0\le p\le q$, if $u r P_{2 i}(u)>\bigl(P_{k-2}(u)/u-(u^2+1)r\bigr)P_{2 i+1}(u)$, then 
\begin{multline*}
g_p\bigl(P_{2 i+1}(u),P_{2 i+3}(u),P_{2 i+k+1}(u)\bigr)\\
=(r-1)P_{2 i+3}(u)+(q+p)P_{2 i+k+1}(u)-P_{2 i+1}(u)\,. 
\end{multline*}  
Otherwise, 
\begin{multline*}
g_p\bigl(P_{2 i+1}(u),P_{2 i+3}(u),P_{2 i+k+1}(u)\bigr)\\
=\left(\frac{P_k(u)}{u}-1\right)P_{2 i+3}(u)+(q+p-1)P_{2 i+k+1}(u)-P_{2 i+1}(u)\,. 
\end{multline*} 
\label{th:i-odd-k-even}
\end{theorem}

\subsection{Examples} 

We shall use the same examples as in Subsection \ref{examples4}. Let $u=2$.  

When $A=\{P_9,P_{11},P_{12}\}$, by $q=197$ and $r=0$, we have for $0\le p\le 99$
$$
g_p(A)=(P_3-1)P_{11}+(2 p+196)P_{12}-P_9\,. 
$$ 
In fact, when $p=100$, this formula yields $5510539$, which does not match the real $g_{100}(A)=5482819$. 

When $A=\{P_7,P_9,P_{11}\}$, by $q=28$ and $r=1$, we have for $0\le p\le 28$
$$
g_p(A)=(P_4/2-1)P_9+(p+27)P_{11}-P_7\,. 
$$ 
In fact, when $p=29$, this formula yields $326252$, which does not match the real $g_{29}(A)=322312$.

\section{$p$-genus for $\bigl(P_{2 i+1}(u),P_{2 i+3}(u),P_{2i+k+1}(u)\bigr)$}  

Let $k$ be odd.  
With reference to Table \ref{tb:ap0podd}, the sum of all of the elements of ${\rm Ap}_p(A)$ is obtained as follows: where it is simplified using the formulas (\ref{eq:rec-k-even}), (\ref{eq:2i+3}) and (\ref{eq:2i+k+1}).
\begin{align*}  
&\sum_{w\in{\rm Ap}_p(A)}w\\
&=\sum_{z=0}^{u-1}\sum_{m=0}^{\mathfrak r-1}\sum_{l=0}^p\left(\bigl(l P_k(u)+m\bigr)P_{2 i+3}(u)+\bigl(\mathfrak q+u(p-2 l)+z\bigr)P_{2 i+k+1}(u)\right)\\
&\quad +\sum_{z=0}^{u-1}\sum_{m=\mathfrak r}^{P_k(u)-1}\sum_{l=0}^{p-1}\left(\bigl(l P_k(u)+m\bigr)P_{2 i+3}(u)\right.\\
&\qquad\qquad\left.+\bigl(\mathfrak q+u(p-2 l-1)+z\bigr)P_{2 i+k+1}(u)\right)\\
&\quad +\sum_{m=0}^{P_k(u)-1}\sum_{l=0}^{\mathfrak q-p u-1}\left(\bigl(p P_k(u)+m\bigr)P_{2 i+3}(u)+l P_{2 i+k+1}(u)\right)\\
&=-\frac{p^2}{2}u P_{2 i+1}(u)P_k(u)P_{k-2}(u)+\frac{p}{2}P_{2 i+1}(u)P_k(u)\bigl(2 P_{2 i+3}(u)-u P_{k-2}(u)\bigr)\\
&\quad +\frac{P_{2 i+1}(u)}{2 u}\biggl(\bigl(P_{2 i+1}(u)-u\bigr)P_{2 i+3}(u)+u(u-1)P_{2 i+k+1}(u)\\
&\qquad -\mathfrak q P_{k-2}(u)\bigl(2 P_{2 i+1}-(\mathfrak q+u)P_{k}(u)\bigr)\biggr)\,. 
\end{align*}
Here, we used the relations (\ref{eq:rec-k-odd}) and 
\begin{equation} 
P_{2 i+3}(u)P_k(u)-u P_{2 i+k+1}(u)=P_{2 i+1}(u)P_{k-2}(u)\,. 
\label{eq:334}
\end{equation}   
Hence, by Lemma \ref{lem-mp} (\ref{mp-n}), we have 
\begin{align*}  
&n_p\bigl(P_{2 i+1}(u),P_{2 i+3}(u),P_{2i+k+1}(u)\bigr)=\frac{1}{P_{2 i+1}(u)}\sum_{w\in{\rm Ap}_p(A)}w-\frac{P_{2 i+1}(u)-1}{2}\\
&=-\frac{p^2}{2}u P_k(u)P_{k-2}(u)+\frac{p}{2}P_k(u)\bigl(2 P_{2 i+3}(u)-u P_{k-2}(u)\bigr)\\
&\quad +\frac{P_{2 i+1}(u)}{2 u}\biggl(\bigl(P_{2 i+1}(u)-u\bigr)\bigl(P_{2 i+3}(u)-u\bigr)+u(u-1)\bigl(P_{2 i+k+1}(u)-1\bigr)\\
&\qquad -\mathfrak q P_{k-2}(u)\bigl(2 P_{2 i+1}-(\mathfrak q+u)P_{k}(u)\bigr)\biggr)\,. 
\end{align*} 

This situation is continued as long as $z=\mathfrak q-p u-1\ge 0$. However, when $p=(\mathfrak q+1)/u$, the shape of the block on the right side collapses. Thus, the regularity of taking the sum of ${\rm Ap}_p(A)$ is broken.  

\begin{theorem}  
Let $i\ge 2$ and $k$ be odd with $3\le k\le 2 i-1$.  
Let $\mathfrak q$ be the largest integer not exceeding $P_{2 i+1}(u)/P_k(u)$ with $\mathfrak q\equiv 1\pmod u$. Then for $0\le p\le(\mathfrak q-1)/u$, 
\begin{align*}  
&n_p\bigl(P_{2 i+1}(u),P_{2 i+3}(u),P_{2i+k+1}(u)\bigr)\\
&=-\frac{p^2}{2}u P_k(u)P_{k-2}(u)+\frac{p}{2}P_k(u)\bigl(2 P_{2 i+3}(u)-u P_{k-2}(u)\bigr)\\
&\quad +\frac{P_{2 i+1}(u)}{2 u}\biggl(\bigl(P_{2 i+1}(u)-u\bigr)\bigl(P_{2 i+3}(u)-u\bigr)+u(u-1)\bigl(P_{2 i+k+1}(u)-1\bigr)\\
&\qquad -\mathfrak q P_{k-2}(u)\bigl(2 P_{2 i+1}-(\mathfrak q+u)P_{k}(u)\bigr)\biggr)\,.  
\end{align*} 
\label{thgenus:i-k-odd}
\end{theorem} 
\bigskip

Let $k$ be even.  
With reference to Table \ref{tb:ap0peven}, the sum of all of the elements of ${\rm Ap}_p(A)$ is obtained as follows: where it is simplified using the formulas (\ref{eq:rec-k-odd}), (\ref{eq:2i+3}) and (\ref{eq:2i+k+1}).
\begin{align*}  
&\sum_{w\in{\rm Ap}_p(A)}w\\
&=\sum_{m=0}^{r-1}\sum_{l=0}^p\left(\left(\frac{l P_k(u)}{u}+m\right)P_{2 i+3}(u)+(q+p-2 l)P_{2 i+k+1}(u)\right)\\
&\quad +\sum_{m=r}^{P_k(u)/u-1}\sum_{l=0}^{p-1}\left(\left(\frac{l P_k(u)}{u}+m\right)P_{2 i+3}(u)+(q+p-2 l-1)P_{2 i+k+1}(u)\right)\\
&\quad +\sum_{m=0}^{P_k(u)/u-1}\sum_{l=0}^{q-p-1}\left(\left(\frac{p P_k(u)}{u}+m\right)P_{2 i+3}(u)+l P_{2 i+k+1}(u)\right)\\
&=-\frac{p^2 P_{2 i+1}(u)P_k(u)P_{k-2}(u)}{2 u^2}+\frac{p P_{2 i+1}(u)}{2 u^2}P_k(u)\bigl(2 u P_{2 i+3}(u)-P_{k-2}(u)\bigr)\\
&\quad +\frac{P_{2 i+1}(u)}{2 u^2}\bigl(u^2\bigl(P_{2 i+1}(u)-1\bigr) P_{2 i+3}(u)\\
&\qquad -q P_{k-2}(u)\bigl(2 u P_{2 i+1}(u)-(q+1)P_k(u)\bigr)\bigr)\,. 
\end{align*} 
Here, we used the relations (\ref{eq:rec-k-even}) and (\ref{eq:334}). 
Hence, by Lemma \ref{lem-mp} (\ref{mp-n}), we have 
\begin{align*}  
&n_p\bigl(P_{2 i+1}(u),P_{2 i+3}(u),P_{2i+k+1}(u)\bigr)=\frac{1}{P_{2 i+1}(u)}\sum_{w\in{\rm Ap}_p(A)}w-\frac{P_{2 i+1}(u)-1}{2}\\
&=-\frac{p^2 P_k(u)P_{k-2}(u)}{2 u^2}+\frac{p}{2 u^2}P_k(u)\bigl(2 u P_{2 i+3}(u)-P_{k-2}(u)\bigr)\\
&\quad +\frac{1}{2 u^2}\bigl(u^2(P_{2 i+1}(u)-1)\bigl(P_{2 i+3}(u)-1\bigr)\\
&\qquad -q P_{k-2}(u)\bigl(2 u P_{2 i+1}(u)-(q+1)P_k(u)\bigr)\bigr)\,. 
\end{align*} 

This situation is regularly continued until $q-p-1\ge 0$, as seen in the right-most block in Table \ref{tb:ap0peven}. When $p=q$, only the rightmost block is shortened to the side, but the other parts appear in accordance with the rules because there is only one (or zero) line as a remainder block. Therefore, the regularity of taking the sum of the elements of ${\rm Ap}_p(A)$ is still maintained. However, when $p=q+1$, the element coming to the bottom left becomes shorter, and the regularity of taking the number is broken.   

\begin{theorem} 
Let $i\ge 2$ and $k$ be even with $4\le k\le 2 i$. \\   
Let $q=\fl{u P_{2 i+1}(u)/P_k(u)}$. Then for $0\le p\le q$ 
\begin{align*}  
&n_p\bigl(P_{2 i+1}(u),P_{2 i+3}(u),P_{2i+k+1}(u)\bigr)\\
&=-\frac{p^2 P_k(u)P_{k-2}(u)}{2 u^2}+\frac{p}{2 u^2}P_k(u)\bigl(2 u P_{2 i+3}(u)-P_{k-2}(u)\bigr)\\
&\quad +\frac{1}{2 u^2}\bigl(u^2(P_{2 i+1}(u)-1)\bigl(P_{2 i+3}(u)-1\bigr)\\
&\qquad -q P_{k-2}(u)\bigl(2 u P_{2 i+1}(u)-(q+1)P_k(u)\bigr)\bigr)\,.  
\end{align*} 
\label{thgenus:i-odd-k-even}
\end{theorem}

\subsection{Examples}  

Let $u=2$.  
For $A:=\{P_7,P_9,P_{10}\}$, by $q=\fl{P_7/P_{3}}=33$ (odd), Theorem \ref{thgenus:i-k-odd} holds for $0\le p\le(33-1)/2=16$.  Indeed, $n_{16}=118324$. But when $p=17$, the formula yields $123079$, which does not match the fact $n_{17}=121704$. 

For $A:=\{P_7,P_9,P_{11}\}$, by $q=\fl{2 P_7/P_{4}}=28$, Theorem \ref{thgenus:i-odd-k-even} holds for $0\le p\le 28$.  Indeed, $n_{28}=243404$. But when $p=29$, the formula yields $249140$, which does not match the fact $n_{29}=244501$.

\section*{Acknowledgement}  

We thank the anonymous referee for careful reading of this manuscript and constructive suggestions.

\end{document}